\documentclass[10pt]{amsart}
\usepackage[T1]{fontenc}
\usepackage{geometry}
\usepackage[latin1] {inputenc}
\usepackage{amsmath}
\usepackage{amsfonts, amssymb, textcomp}
\usepackage[colorlinks=flase, linkcolor=red,urlcolor=green, citecolor=blue]{hyperref}
\usepackage{subeqnarray}
\usepackage{color}

\usepackage{latexsym}
\usepackage{fancyhdr}
\usepackage{longtable}
\usepackage{amsmath, amssymb}
\usepackage{graphicx}

\numberwithin{equation}{section}

\newtheorem{lemma}[subsection]{Lemma}
\newtheorem{theorem}[subsection]{Theorem}
\newtheorem{proposition}[subsection]{Proposition}
\newtheorem{corollary}[subsection]{Corollary}
\newtheorem{definition}[subsection]{Definition}

\newtheorem{remark}[subsection]{Remark}

\newcommand{\RR}{\mathbb{R}}
\newcommand{\CC}{\mathbb{C}}
\newcommand{\NN}{\mathbb{N}}

\newcommand{\supp}{\operatorname{supp}}

\let\on=\operatorname

\newenvironment{demo}[1]{\par\smallskip\noindent{\bf #1.}}{\par\smallskip}

\title[Characterization of ultradifferentiable test functions]
{Characterization of ultradifferentiable test functions defined by weight matrices in terms of their Fourier Transform}

\author[G.~Schindl]{Gerhard Schindl}

\address{G.~Schindl: Fakult\"at f\"ur Mathematik, Universit\"at Wien,
Oskar-Morgenstern-Platz~1, A-1090 Wien, Austria}
\email{a0304518@unet.univie.ac.at}

\begin{document}

\begin{abstract}
We prove that functions with compact support in non-quasianalytic classes $\mathcal{E}_{\{\mathcal{M}\}}$ of Roumieu-type and $\mathcal{E}_{(\mathcal{M})}$ of Beurling-type defined by a weight matrix $\mathcal{M}$ with some mild regularity conditions can be characterized by the decay properties of their Fourier transform. For this we introduce the abstract technique of constructing from $\mathcal{M}$ multi-index matrices and associated function spaces. We study the behaviour of this construction in detail and characterize its stability. Moreover non-quasianalyticity of the classes $\mathcal{E}_{\{\mathcal{M}\}}$ and $\mathcal{E}_{(\mathcal{M})}$ is characterized.
\end{abstract}

\thanks{GS was supported by FWF-Project P~23028-N13 and FWF-Project P~26735-N25}
\date{\today}

\maketitle
\section{Introduction}
Spaces of ultradifferentiable functions are sub-classes of smooth functions with certain growth conditions on all their derivatives. In the literature two different approaches are considered to introduce these classes, either using a weight sequence $M=(M_k)_k$ or using a weight function $\omega$. Given a compact set $K$ the classes
$$\left\{\frac{f^{(k)}(x)}{h^kM_k}: x\in K, k\in\NN\right\}\hspace{30pt}\text{respectively}\hspace{30pt}\left\{\frac{f^{(k)}(x)}{\exp(1/l\varphi^{*}_{\omega}(lk))}: x\in K, k\in\NN\right\}$$
should be bounded, where the positive real number $h$ or $l$ is subject to either a universal or an existential quantifier and $\varphi^{*}_{\omega}$ denotes the Young-conjugate of $\varphi_{\omega}=\omega\circ\exp$. In the case of a universal quantifier we call the class of {\itshape Beurling type}, denoted by $\mathcal{E}_{(M)}$ or $\mathcal{E}_{(\omega)}$. In the case of an existential quantifier we call the class of {\itshape Roumieu type}, denoted by $\mathcal{E}_{\{M\}}$ or $\mathcal{E}_{\{\omega\}}$. In the following we write $\mathcal{E}_{[\star]}$ if either $\mathcal{E}_{\{\star\}}$ or $\mathcal{E}_{(\star)}$ is considered.\vspace{6pt}

The classes $\mathcal{E}_{[M]}$ were considered earlier than $\mathcal{E}_{[\omega]}$. For the weight sequence approach see e.g. \cite{mandelbrojtbook} and \cite{Komatsu73}, for $\mathcal{E}_{[\omega]}$ we refer to \cite{BraunMeiseTaylor90}. In \cite{BonetMeiseMelikhov07} both methods were compared and it was shown that in general a class $\mathcal{E}_{[M]}$ cannot be obtained by a weight function $\omega$ and vice versa. At the beginning, ultradifferentiable classes were
studied using the growth of the derivatives and later with the Fourier transform. Finally,
Braun, Meise and Taylor in \cite{BraunMeiseTaylor90} have unified both theories. For a detailed survey we refer to the introductions in \cite{BraunMeiseTaylor90} and \cite{BonetMeiseMelikhov07}.\vspace{6pt}

In \cite{compositionpaper} we have considered classes $\mathcal{E}_{[\mathcal{M}]}$ defined by (one-parameter) weight matrices $\mathcal{M}:=\{M^x: x\in\Lambda\}$. The spaces $\mathcal{E}_{[M]}$ and $\mathcal{E}_{[\omega]}$ were identified as particular cases of $\mathcal{E}_{[\mathcal{M}]}$ but one is able to describe more classes, e.g. the class defined by the {\itshape Gevrey-matrix} $\mathcal{G}:=\{(p!^{s+1})_{p\in\NN}: s>0\}$, see \cite[5.19]{compositionpaper}. Using this new method one is able to transfer results from one setting into the other one and to prove results for $\mathcal{E}_{[M]}$ and $\mathcal{E}_{[\omega]}$ simultaneously, e.g. see \cite{compositionpaper} and \cite{characterizationstabilitypaper}.\vspace{6pt}

The main aim of this work is to show that assuming some mild properties for $\mathcal{M}$ the functions with compact support $\mathcal{D}_{[\mathcal{M}]}\subseteq\mathcal{E}_{[\mathcal{M}]}$ can be characterized in terms of the decay properties of their Fourier transform.\vspace{6pt}

First, we generalize in Section \ref{section3} a central new idea in \cite{compositionpaper}. We have shown that to each $\omega$ we can associate a weight matrix $\Omega:=\{(\Omega^l_j)_{j\ge 0}: l>0\}$, defined by $\Omega^l_j:=\exp(1/l\varphi^{*}_{\omega}(lj))$, such that $\mathcal{E}_{[\omega]}=\mathcal{E}_{[\Omega]}$ holds as locally convex vector spaces.

In this work we start with an abstractly given weight matrix $\mathcal{M}=\{M^x: x\in\Lambda\}$ satisfying some standard assumptions. To $\mathcal{M}$ we associate another matrix $\omega_{\mathcal{M}}:=\{\omega_{M^x}: x\in\Lambda\}$ consisting of {\itshape associated functions} $\omega_{M^x}$. Applying again the idea of \cite{compositionpaper} we obtain a matrix $\{M^{x;l}: x\in\Lambda, l>0\}$ and iterating this procedure we get a sequence of multi-index weight matrices consisting of weight sequences and weight functions. In Section \ref{section3} this technique is studied in detail.\vspace{6pt}

First, in \ref{subsection31}, we will characterize the case where all multi-index weight matrices of weight sequences are equivalent. Thus $\mathcal{E}_{[\mathcal{M}]}$ is stable as locally convex vector space under adjoining indices, see Theorem \ref{stablefirstcentralresult}. It will turn out that only in the first step a non-stable effect can occur, see Corollary \ref{1consequ}.

The spaces associated to the matrices of weight functions in this construction are always stable. Using results from \ref{subsection32} and Theorem \ref{stablefirstcentralresult} we can prove the first main result Theorem \ref{4}: As locally convex vector spaces the equality $\mathcal{E}_{[\mathcal{M}]}=\mathcal{E}_{[\omega_{\mathcal{M}}]}$ is valid.\vspace{6pt}

In the next step, in Section \ref{section4}, we characterize the non-quasianalyticity of $\mathcal{E}_{[\mathcal{M}]}$, see Theorem \ref{Matrix-non-quasi-analyticity}. Thus the cases where the spaces $\mathcal{D}_{[\mathcal{M}]}$ are non-trivial are classified. The Roumieu case is quite clear and for the Beurling case we generalize \cite[Lemma 5.1]{intersectionpaper}, where stronger conditions for the matrix $\mathcal{M}$ were assumed.\vspace{6pt}

In Section \ref{section5} we combine Theorem \ref{4} and Theorem \ref{Matrix-non-quasi-analyticity}. Using and generalizing the methods and estimates introduced in \cite{BraunMeiseTaylor90} we are able to characterize functions in $\mathcal{D}_{[\mathcal{M}]}$ in terms of the decay properties of their Fourier transform, see Theorem \ref{centralfouriertheorem}. As special case this holds for the Gevrey-matrix $\mathcal{G}$.\vspace{6pt}

Finally, in Section \ref{section6}, we apply the technique of associating a weight matrix to prove some variations of comparison results due to \cite{BonetMeiseMelikhov07} concerning the classes $\mathcal{E}_{[M]}$ and $\mathcal{E}_{[\omega]}$.\vspace{6pt}

This work contains some results of the author PhD Thesis, see \cite{dissertation}. The author thanks his advisors A. Kriegl, P.W. Michor and A. Rainer for the supervision and their helpful ideas. 

\subsection{Basic notation}
We denote by $\mathcal{E}$ the class of smooth functions, $\mathcal{C}^{\omega}$ is the class of all real analytic functions. We will write $\NN_{>0}=\{1,2,\dots\}$ and $\NN=\NN_{>0}\cup\{0\}$. Moreover we put $\RR_{>0}:=\{x\in\RR: x>0\}$, i.e. the set of all positive real numbers. For $\alpha=(\alpha_1,\dots,\alpha_n)\in\NN^n$ we use the usual multi-index notation, write $\alpha!:=\alpha_1!\dots\alpha_n!$, $|\alpha|:=\alpha_1+\dots+\alpha_n$ and for $x=(x_1,\dots,x_n)\in\RR^n$ we set $x^{\alpha}=x_1^{\alpha_1}\cdots x_n^{\alpha_n}$. We also put $\partial^{\alpha}=\partial_1^{\alpha_1}\cdots\partial_n^{\alpha_n}$ and for a given function $f:U\subseteq\RR^r\rightarrow\RR^s$ defined on a non-empty open set $U\subseteq\RR^r$ we denote by $f^{(k)}$ the $k$-th order {\itshape Fréchet derivative} of $f$. Let $E_1,\dots,E_k$ and $F$ be topological vector spaces, then $L(E_1,\dots,E_k,F)$ is the space of all bounded $k$-linear mappings $E_1\times\dots\times E_k\rightarrow F$. If $E=E_i$ for $i=1,\dots,k$, then we write $L^k(E,F)$. With $\|\cdot\|_{\RR^n}$ we denote the Euclidian norm on $\RR^n$.

Let $K\subset\subset\RR^r$ be a compact set with smooth boundary, then $\mathcal{E}(K,\RR^s)$ denotes the space of all smooth functions on the interior $K^{\circ}$ such that each derivative of $f$ can be continuously extended to $K$.

{\itshape Convention:} Let $\star\in\{M,\omega,\mathcal{M}\}$, then we write $\mathcal{E}_{[\star]}$ if either $\mathcal{E}_{\{\star\}}$ or $\mathcal{E}_{(\star)}$ is considered with the following restriction: Statements that involve more than one $\mathcal{E}_{[\star]}$ symbol must not be interpreted by mixing $\mathcal{E}_{\{\star\}}$ and $\mathcal{E}_{(\star)}$. The same notation resp. convention will be used for the conditions, so write $(\mathcal{M}_{[\star]})$ for either $(\mathcal{M}_{\{\star\}})$ or $(\mathcal{M}_{(\star)})$.

\section{Basic definitions}
\subsection{Weight sequences and classes of ultradifferentiable functions $\mathcal{E}_{[M]}$}
$M=(M_k)_k\in\RR_{>0}^{\NN}$ is called a {\itshape weight sequence}. We introduce also $m=(m_k)_k$ defined by $m_k:=\frac{M_k}{k!}$ and $\mu=(\mu_k)_k$ by $\mu_k:=\frac{M_k}{M_{k-1}}$, $\mu_0:=1$. $M$ is called {\itshape normalized} if $1=M_0\le M_1$ holds (w.l.o.g.).

$(1)$ $M$ is called {\itshape log-convex} if
$$\hypertarget{lc}{(\text{lc})}:\Leftrightarrow\;\forall\;j\in\NN:\;M_j^2\le M_{j-1} M_{j+1}.$$
$M$ is log-convex if and only if $(\mu_k)_k$ is increasing. If $M$ is log-convex and normalized, then $M$ and $k\mapsto(M_k)^{1/k}$ are both increasing, see e.g. \cite[Lemma 2.0.4]{diploma}.

$(2)$ $M$ has {\itshape moderate growth} if
$$\hypertarget{mg}{(\text{mg})}:\Leftrightarrow\exists\;C\ge 1\;\forall\;j,k\in\NN:\;M_{j+k}\le C^{j+k}M_j M_k.$$

$(3)$ $M$ is called {\itshape non-quasianalytic} if
$$\hypertarget{mnq}{(\text{nq})}:\Leftrightarrow\;\sum_{p=1}^{\infty}\frac{M_{p-1}}{M_p}<+\infty.$$
Using {\itshape Carleman's inequality} one can show that if $M$ has \hyperlink{lc}{$(\text{lc})$}, then $$\sum_{p=1}^{\infty}\frac{M_{p-1}}{M_p}<+\infty\Leftrightarrow\sum_{p=1}^{\infty}\frac{1}{(M_p)^{1/p}}<+\infty.$$

$(4)$ $M$ has $(\beta_3)$ if
$$\exists\;Q\in\NN_{>0}:\;\liminf_{p\rightarrow\infty}\frac{\mu_{Qp}}{\mu_p}>1.$$

$(5)$ For $M=(M_p)_p$ and $N=(N_p)_p$ we write $M\le N$ if and only if $M_p\le N_p$ holds for all $p\in\NN$. Moreover we define
$$M\hypertarget{mpreceq}{\preceq}N:\Leftrightarrow\;\exists\;C_1,C_2\ge 1\;\forall\;p\in\NN:\; M_p\le C_2 C_1^p N_p\Longleftrightarrow\sup_{p\in\NN_{>0}}\left(\frac{M_p}{N_p}\right)^{1/p}<+\infty$$
and call the sequences equivalent if
$$M\hypertarget{approx}{\approx}N:\Leftrightarrow\;M\hyperlink{mpreceq}{\preceq}N\;\text{and}\;N\hyperlink{mpreceq}{\preceq}M.$$
\hyperlink{mg}{$(\text{mg})$} and \hyperlink{mnq}{$(\text{nq})$} are stable w.r.t. \hyperlink{approx}{$\approx$}. Furthermore we will write
$$M\hypertarget{mtriangle}{\vartriangleleft}N:\Leftrightarrow\;\forall\;h>0\;\exists\;C_h\ge 1\;\forall\;j\in\NN:\; M_j\le C_h h^j N_j\Longleftrightarrow\lim_{p\rightarrow\infty}\left(\frac{M_p}{N_p}\right)^{1/p}=0.$$
For convenience we introduce the set
$$\hypertarget{LCset}{\mathcal{LC}}:=\{M\in\RR_{>0}^{\NN}:\;M\;\text{is normalized, log-convex},\;\lim_{k\rightarrow\infty}(M_k)^{1/k}=+\infty\}.$$
Let $r,s\in\NN_{>0}$ and $U\subseteq\RR^r$ be a non-empty open set. We introduce the classes of ultradifferentiable functions of Roumieu type by
$$\mathcal{E}_{\{M\}}(U,\RR^s):=\{f\in\mathcal{E}(U,\RR^s):\;\forall\;K\subset\subset U\;\exists\;h>0:\;\|f\|_{M,K,h}<+\infty\},$$
and the classes of ultradifferentiable functions of Beurling type by
$$\mathcal{E}_{(M)}(U,\RR^s):=\{f\in\mathcal{E}(U,\RR^s):\;\forall\;K\subset\subset U\;\forall\;h>0:\;\|f\|_{M,K,h}<+\infty\},$$
where we denote
\begin{equation}\label{semi-norm-2}
\|f\|_{M,K,h}:=\sup_{k\in\NN,x\in K}\frac{\|f^{(k)}(x)\|_{L^k(\RR^r,\RR^s)}}{h^{k} M_{k}}
\end{equation}
and $\|f^{(k)}(x)\|_{L^k(\RR^r,\RR^s)}:=\sup\{\|f^{(k)}(x)(v_1,\dots,v_k)\|_{\RR^s}: \|v_i\|_{\RR^r}\le 1\;\forall\;1\le i\le k\}$.

For a compact set $K$ with smooth boundary
$$\mathcal{E}_{M,h}(K,\RR^s):=\{f\in\mathcal{E}(K,\RR^s): \|f\|_{M,K,h}<+\infty\}$$
is a Banach space and we define the following topological vector spaces
\begin{equation}\label{repr3}
\mathcal{E}_{\{M\}}(U,\RR^s):=\underset{K\subset\subset U}{\varprojlim}\;\underset{h>0}{\varinjlim}\;\mathcal{E}_{M,h}(K,\RR^s)=\underset{K\subseteq U}{\varprojlim}\;\mathcal{E}_{\{M\}}(K,\RR^s)
\end{equation}
and
\begin{equation}\label{repr4}
\mathcal{E}_{(M)}(U,\RR^s):=\underset{K\subset\subset U}{\varprojlim}\;\underset{h>0}{\varprojlim}\;\mathcal{E}_{M,h}(K,\RR^s)=\underset{K\subseteq U}{\varprojlim}\;\mathcal{E}_{(M)}(K,\RR^s).
\end{equation}
In $\mathcal{E}_{M,h}(K,\RR^s)$ instead of compact sets $K$ with smooth boundary one can also consider a relatively compact open subset $K$ of $U$ (see \cite{thilliez}) or one can work with {\itshape Whitney jets} on the compact set $K$ (see \cite{Komatsu73} and also \cite{BonetMeiseMelikhov07}).

We recall some facts for log-convex $M$:

\begin{itemize}
\item[$(i)$] We write $\mathcal{E}^{\text{global}}_{\{M\}}(U,\RR^s):=\{f\in\mathcal{E}(U,\RR^s):\;\exists\;h>0:\;\|f\|_{M,U,h}<+\infty\}$. Then there exist {\itshape characteristic functions}
\begin{equation}\label{characteristicfunction}
\theta_M\in\mathcal{E}^{\text{global}}_{\{M\}}(\RR,\RR): \forall\;j\in\NN: \left|\theta_M^{(j)}(0)\right|\ge M_j,
\end{equation}
see \cite[Lemma 2.9]{compositionpaper} and \cite[Theorem 1]{thilliez}. Note that the Beurling class $\mathcal{E}^{\text{global}}_{(M)}(\RR,\RR)$ cannot contain such $\theta_M$, see \cite[Proposition 3.1.2]{diploma}.

\item[$(ii)$] If $N$ is arbitrary, then $M\hyperlink{mpreceq}{\preceq} N\Longleftrightarrow\mathcal{E}_{\{M\}}\subseteq\mathcal{E}_{\{N\}}$ and $M\hyperlink{mtriangle}{\vartriangleleft}N\Longleftrightarrow\mathcal{E}_{\{M\}}\subseteq\mathcal{E}_{(N)}$. If $M\in\hyperlink{LCset}{\mathcal{LC}}$, then $M\hyperlink{mpreceq}{\preceq} N\Longleftrightarrow\mathcal{E}_{[M]}\subseteq\mathcal{E}_{[N]}$.

\item[$(iii)$] For any non-empty open set $U\subseteq\RR^r$ both classes $\mathcal{E}_{\{M\}}(U,\RR)$ and $\mathcal{E}_{(M)}(U,\RR)$ are closed under pointwise multiplication, see e.g. \cite[Proposition 2.0.8]{diploma}.
\end{itemize}

\subsection{Classes of ultradifferentiable functions defined by weight matrices}\label{weighmatrixdefinition}
\begin{definition}
Let $(\Lambda,\le)$ be a partially ordered set which is both up- and downward directed, $\Lambda=\RR_{>0}$ is the most important example. A {\itshape weight matrix} $\mathcal{M}$ associated to $\Lambda$ is a family of weight sequences $\mathcal{M}:=\{M^x\in\RR_{>0}^{\NN}: x\in\Lambda\}$ such that
$$\hypertarget{Marb}{(\mathcal{M})}:\Leftrightarrow\;\forall\;x\in\Lambda:\;M^x\;\text{is normalized, increasing},\;M^x\le M^y\;\text{for}\;x\le y.$$
We call $\mathcal{M}$ {\itshape standard log-convex,} if
$$\hypertarget{Msc}{(\mathcal{M}_{\on{sc}})}:\Leftrightarrow(\mathcal{M})\;\text{and}\;\forall\;x\in\Lambda:\;M^x\in\hyperlink{LCset}{\mathcal{LC}}.$$
Also the sequences $m^x_k:=\frac{M^x_k}{k!}$ and $\mu^x_k:=\frac{M^x_k}{M^x_{k-1}}$, $\mu^x_0:=1$, will be used.
\end{definition}

We introduce spaces of vector-valued ultradifferentiable functions classes defined by a weight matrices of Roumieu type $\mathcal{E}_{\{\mathcal{M}\}}$ and Beurling type $\mathcal{E}_{(\mathcal{M})}$ as follows, see also \cite[4.2]{compositionpaper}.

Let $r,s\in\NN_{>0}$, let $U\subseteq\RR^r$ be a non-empty open set. For all compact sets $K\subset\subset U$ we put
\begin{equation}\label{generalroumieu}
\mathcal{E}_{\{\mathcal{M}\}}(K,\RR^s):=\bigcup_{x\in\Lambda}\mathcal{E}_{\{M^x\}}(K,\RR^s)\hspace{20pt}\mathcal{E}_{\{\mathcal{M}\}}(U,\RR^s):=\bigcap_{K\subset\subset U}\bigcup_{x\in\Lambda}\mathcal{E}_{\{M^x\}}(K,\RR^s)
\end{equation}
and
\begin{equation}\label{generalbeurling}
\mathcal{E}_{(\mathcal{M})}(K,\RR^s):=\bigcap_{x\in\Lambda}\mathcal{E}_{(M^x)}(K,\RR^s)\hspace{20pt}\mathcal{E}_{(\mathcal{M})}(U,\RR^s):=\bigcap_{x\in\Lambda}\mathcal{E}_{(M^x)}(U,\RR^s).
\end{equation}
For a compact set $K\subset\subset\RR^r$ one has the representations
$$\mathcal{E}_{\{\mathcal{M}\}}(K,\RR^s):=\underset{x\in\Lambda}{\varinjlim}\;\underset{h>0}{\varinjlim}\;\mathcal{E}_{M^x,h}(K,\RR^s)$$
and so for $U\subseteq\RR^r$ non-empty open
\begin{equation}\label{generalroumieu1}
\mathcal{E}_{\{\mathcal{M}\}}(U,\RR^s):=\underset{K\subset\subset U}{\varprojlim}\;\underset{x\in\Lambda}{\varinjlim}\;\underset{h>0}{\varinjlim}\;\mathcal{E}_{M^x,h}(K,\RR^s).
\end{equation}
Similarly we get for the Beurling case
\begin{equation}\label{generalbeurling1}
\mathcal{E}_{(\mathcal{M})}(U,\RR^s):=\underset{K\subset\subset U}{\varprojlim}\;\underset{x\in\Lambda}{\varprojlim}\;\underset{h>0}{\varprojlim}\;\mathcal{E}_{M^x,h}(K,\RR^s).
\end{equation}

If $\Lambda=\RR_{>0}$ we can assume that all occurring limits are countable and restrict to $\Lambda=\NN_{>0}$ in the Roumieu case. Thus $\mathcal{E}_{(\mathcal{M})}(U,\RR^s)$ is a {\itshape Fréchet space} and $\underset{x\in\Lambda}{\varinjlim}\;\underset{h>0}{\varinjlim}\;\mathcal{E}_{M^x,h}(K,\RR^s)=\underset{n\in\NN_{>0}}{\varinjlim}\;\;\mathcal{E}_{M^n,n}(K,\RR^s)$ is a {\itshape Silva space}, i.e. a countable inductive limit of Banach spaces with compact connecting mappings. For more details concerning the locally convex topology we refer to \cite[4.2-4.4]{compositionpaper}. In the appendix in Proposition \ref{Nuclearproposition} we will show that for some weight matrices the connecting mappings are even nuclear.

\subsection{Conditions for a weight matrix $\mathcal{M}$}
We are going to introduce several conditions on $\mathcal{M}$, see also \cite[4.1]{compositionpaper}. First consider the following conditions {\itshape of Roumieu type}.
\begin{itemize}
\item[\hypertarget{R-dc}{$(\mathcal{M}_{\{\text{dc}\}})$}] $\forall\;x\in\Lambda\;\exists\;C>0\;\exists\;y\in\Lambda\;\forall\;j\in\NN: M^x_{j+1}\le C^{j+1} M^y_j$
\item[\hypertarget{R-mg}{$(\mathcal{M}_{\{\text{mg}\}})$}] $\forall\;x\in\Lambda\;\exists\;C>0\;\exists\;y_1,y_2\in\Lambda\;\forall\;j,k\in\NN: M^x_{j+k}\le C^{j+k} M^{y_1}_j M^{y_2}_k$
\item[\hypertarget{R-L}{$(\mathcal{M}_{\{\text{L}\}})$}] $\forall\;C>0\;\forall\;x\in\Lambda\;\exists\;D>0\;\exists\;y\in\Lambda\;\forall\;k\in\NN: C^k M^x_k\le D M^y_k$
\item[\hypertarget{R-strict}{$(\mathcal{M}_{\{\text{strict}\}})$}] $\forall\;x\in\Lambda\;\exists\;y\in\Lambda\;:\;\;\sup_{k\in\NN_{>0}}\left(\frac{M^y_k}{M^x_k}\right)^{1/k}=+\infty$
\item[\hypertarget{R-BR}{$(\mathcal{M}_{\{\text{BR}\}})$}] $\forall\;x\in\Lambda\;\exists\;y\in\Lambda: M^x\hyperlink{mtriangle}{\vartriangleleft}M^y$
\end{itemize}

Analogously we introduce the {\itshape Beurling type conditions}.
\begin{itemize}
\item[\hypertarget{B-dc}{$(\mathcal{M}_{(\text{dc})})$}] $\forall\;x\in\Lambda\;\exists\;C>0\;\exists\;y\in\Lambda\;\forall\;j\in\NN: M^y_{j+1}\le C^{j+1} M^x_j$
\item[\hypertarget{B-mg}{$(\mathcal{M}_{(\text{mg})})$}] $\forall\;x_1,x_2\in\Lambda\;\exists\;C>0\;\exists\;y\in\Lambda\;\forall\;j,k\in\NN: M^y_{j+k}\le C^{j+k} M^{x_1}_j M^{x_2}_k$
\item[\hypertarget{B-L}{$(\mathcal{M}_{(\text{L})})$}] $\forall\;C>0\;\forall\;x\in\Lambda\;\exists\;D>0\;\exists\;y\in\Lambda\;\forall\;k\in\NN: C^k M^y_k\le D M^x_k$
\item[\hypertarget{B-strict}{$(\mathcal{M}_{(\text{strict})})$}] $\forall\;x\in\Lambda\;\exists\;y\in\Lambda\;:\;\;\sup_{k\in\NN_{>0}}\left(\frac{M^x_k}{M^y_k}\right)^{1/k}=+\infty$
\item[\hypertarget{B-BR}{$(\mathcal{M}_{(\text{BR})})$}] $\forall\;x\in\Lambda\;\exists\;y\in\Lambda: M^y\hyperlink{mtriangle}{\vartriangleleft}M^x$
\end{itemize}

\subsection{Inclusion relations}
Given two matrices $\mathcal{M}=\{M^x: x\in\Lambda\}$ and $\mathcal{N}=\{N^y: y\in\Lambda'\}$ we introduce
$$\mathcal{M}\hypertarget{Mroumpreceq}{\{\preceq\}}\mathcal{N}\;:\Leftrightarrow\forall\;x\in\Lambda\;\exists\;y\in\Lambda':\;M^x\hyperlink{mpreceq}{\preceq}N^y$$
and
$$\mathcal{M}\hypertarget{Mbeurpreceq}{(\preceq)}\mathcal{N}\;:\Leftrightarrow\forall\;y\in\Lambda'\;\exists\;x\in\Lambda:\;M^x\hyperlink{mpreceq}{\preceq}N^y.$$
By definition $\mathcal{M}[\preceq]\mathcal{N}$ implies $\mathcal{E}_{[\mathcal{M}]}\subseteq\mathcal{E}_{[\mathcal{N}]}$ and write
$$\mathcal{M}\hypertarget{Mroumapprox}{\{\approx\}}\mathcal{N}\;:\Leftrightarrow\mathcal{M}\hyperlink{Mroumpreceq}{\{\preceq\}}\mathcal{N}\;\text{and}\;\mathcal{N}\hyperlink{Mroumpreceq}{\{\preceq\}}\mathcal{M}$$
and
$$\mathcal{M}\hypertarget{Mbeurapprox}{(\approx)}\mathcal{N}\;:\Leftrightarrow\mathcal{M}\hyperlink{Mbeurpreceq}{(\preceq)}\mathcal{N}\;\text{and}\;\mathcal{N}\hyperlink{Mbeurpreceq}{(\preceq)}\mathcal{M}.$$
Moreover, we introduce
$$\mathcal{M}\hypertarget{Mtriangle}{\vartriangleleft}\mathcal{N}\;:\Leftrightarrow\forall\;x\in\Lambda\;\forall\;y\in\Lambda': M^x\hyperlink{mtriangle}{\vartriangleleft}N^{y},$$
so $\mathcal{M}\hypertarget{Mtriangle}{\vartriangleleft}\mathcal{N}$ implies $\mathcal{E}_{\{\mathcal{M}\}}\subseteq\mathcal{E}_{(\mathcal{N})}$. In \cite[Proposition 4.6]{compositionpaper} the relations above were characterized for \hyperlink{Msc}{$(\mathcal{M}_{\on{sc}})$} matrices with $\Lambda=\Lambda'=\RR_{>0}$. In this context we introduce also

$\hypertarget{R-Comega}{(\mathcal{M}_{\{\mathcal{C}^{\omega}\}})}\hspace{5pt}\exists\;x\in\Lambda:\;\liminf_{k\rightarrow\infty}(m^x_k)^{1/k}>0$,

$\hypertarget{holom}{(\mathcal{M}_{\mathcal{H}})}\hspace{16pt}\forall\;x\in\Lambda:\;\liminf_{k\rightarrow\infty}(m^x_k)^{1/k}>0$,

$\hypertarget{B-Comega}{(\mathcal{M}_{(\mathcal{C}^{\omega})})}\hspace{7pt} \forall\;x\in\Lambda:\;\lim_{k\rightarrow\infty}(m^x_k)^{1/k}=+\infty$.\vspace{6pt}

Recall \cite[Proposition 4.6]{compositionpaper}: If \hyperlink{R-Comega}{$(\mathcal{M}_{\{\mathcal{C}^{\omega}\}})$} holds then the class of real-analytic-functions is contained in $\mathcal{E}_{\{\mathcal{M}\}}$, if \hyperlink{B-Comega}{$(\mathcal{M}_{(\mathcal{C}^{\omega})})$} then the real-analytic functions are contained in $\mathcal{E}_{(\mathcal{M})}$. If \hyperlink{holom}{$(\mathcal{M}_{\mathcal{H}})$} is satisfied, then the restrictions of entire functions are contained in $\mathcal{E}_{(\mathcal{M})}$.

{\itshape Convention:} If $\Lambda=\RR_{>0}$ or $\NN_{>0}$, then $\RR_{>0}$ or $\NN_{>0}$ are always regarded with its natural order $\le$. We will call $\mathcal{M}$ {\itshape constant} if $\mathcal{M}=\{M\}$ or more generally if $M^x\hyperlink{approx}{\approx}M^y$ for all $x,y\in\Lambda$, which violates both \hyperlink{R-strict}{$(\mathcal{M}_{\{\on{strict}\}})$} and \hyperlink{B-strict}{$(\mathcal{M}_{(\on{strict})})$}. Otherwise it will be called {\itshape non-constant}.

\subsection{Classes of ultradifferentiable functions $\mathcal{E}_{[\omega]}$}
A function $\omega:[0,\infty)\rightarrow[0,\infty)$ (sometimes $\omega$ is extended to $\CC$, by $\omega(x):=\omega(|x|)$) is called a {\itshape weight function} if
\begin{itemize}
\item[$(i)$] $\omega$ is continuous,
\item[$(ii)$] $\omega$ is increasing,
\item[$(iii)$] $\omega(x)=0$ for all $x\in[0,1]$ (normalization, w.l.o.g.),
\item[$(iv)$] $\lim_{x\rightarrow\infty}\omega(x)=+\infty$.
\end{itemize}
For convenience we will write that $\omega$ has $\hypertarget{om0}{(\omega_0)}$ if it satisfies $(i)-(iv)$.

Moreover we consider the following conditions:
\begin{itemize}
\item[\hypertarget{om1}{$(\omega_1)}$] $\omega(2t)=O(\omega(t))$ as $t\rightarrow+\infty$.

\item[\hypertarget{om2}{$(\omega_2)$}] $\omega(t)=O(t)$ as $t\rightarrow\infty$.

\item[\hypertarget{om3}{$(\omega_3)$}] $\log(t)=o(\omega(t))$ as $t\rightarrow+\infty$ ($\Leftrightarrow\lim_{t\rightarrow+\infty}\frac{t}{\varphi_{\omega}(t)}=0$).

\item[\hypertarget{om4}{$(\omega_4)$}] $\varphi_{\omega}:t\mapsto\omega(e^t)$ is a convex function on $\RR$.

\item[\hypertarget{om5}{$(\omega_5)$}] $\omega(t)=o(t)$ as $t\rightarrow+\infty$.

\item[\hypertarget{om6}{$(\omega_6)$}] $\exists\;H\ge 1\;\forall\;t\ge 0:\;2\omega(t)\le\omega(H t)+H$.

\item[\hypertarget{om7}{$(\omega_7)$}] $\exists\;H>0\;\exists\;C>0\;\forall\;t\ge 0: \omega(t^2)\le C\omega(H t)+C.$

\item[\hypertarget{omnq}{$(\omega_{\text{nq}})$}] $\int_1^{\infty}\frac{\omega(t)}{t^2}dt<\infty.$
\end{itemize}
An interesting example is $\omega_s(t):=\max\{0,\log(t)^s\}$, $s>1$, which satisfies all listed properties {\itshape except} \hyperlink{om6}{$(\omega_6)$}. For convenience we define the sets
$$\hypertarget{omset0}{\mathcal{W}_0}:=\{\omega:[0,\infty)\rightarrow[0,\infty): \omega\;\text{has}\;\hyperlink{om0}{(\omega_0)},\hyperlink{om3}{(\omega_3)},\hyperlink{om4}{(\omega_4)}\},$$
$$\hypertarget{omset1}{\mathcal{W}}:=\{\omega\in\mathcal{W}_0: \omega\;\text{has}\;\hyperlink{om1}{(\omega_1)}\}.$$
For $\omega\in\hyperlink{omset0}{\mathcal{W}_0}$ we can define the {\itshape Legendre-Fenchel-Young-conjugate} $\varphi^{*}_{\omega}$ by
$$\varphi^{*}_{\omega}(x):=\sup\{x y-\varphi_{\omega}(y): y\ge 0\},$$
with the following properties, e.g. see \cite[Remark 1.3, Lemma 1.5]{BraunMeiseTaylor90}: It is convex and increasing, $\varphi^{*}_{\omega}(0)=0$, $\varphi^{**}_{\omega}=\varphi_{\omega}$, $\lim_{x\rightarrow\infty}\frac{x}{\varphi^{*}_{\omega}(x)}=0$ and finally $x\mapsto\frac{\varphi_{\omega}(x)}{x}$ and $x\mapsto\frac{\varphi^{*}_{\omega}(x)}{x}$ are increasing on $[0,+\infty)$ .

For two weights $\sigma,\tau\in\hyperlink{omset0}{\mathcal{W}_0}$ we write
$$\sigma\hypertarget{ompreceq}{\preceq}\tau:\Leftrightarrow\tau(t)=O(\sigma(t))\;\text{as}\;t\rightarrow+\infty$$
and call them equivalent if
$$\sigma\hypertarget{sim}{\sim}\tau:\Leftrightarrow\sigma\hyperlink{ompreceq}{\preceq}\tau\;\text{and}\;\tau\hyperlink{ompreceq}{\preceq}\sigma.$$
Moreover introduce
$$\sigma\hypertarget{omtriangle}{\vartriangleleft}\tau:\Leftrightarrow\tau(t)=o(\sigma(t))\;\text{as}\;t\rightarrow+\infty.$$
Let $r,s\in\NN_{>0}$, $U\subseteq\RR^r$ be a non-empty open set and $\omega\in\hyperlink{omset0}{\mathcal{W}_0}$. The space of vector-valued ultradifferentiable functions of Roumieu type is defined by
$$\mathcal{E}_{\{\omega\}}(U,\RR^s):=\{f\in\mathcal{E}(U,\RR^s):\;\forall\;K\subset\subset U\;\exists\;l>0:\;\|f\|_{\omega,K,l}<+\infty\}$$
and the space of vector-valued ultradifferentiable functions of Beurling type by
$$\mathcal{E}_{(\omega)}(U,\RR^s):=\{f\in\mathcal{E}(U,\RR^s):\;\forall\;K\subset\subset U\;\forall\;l>0:\;\|f\|_{\omega,K,l}<+\infty\},$$
where
\begin{equation}\label{semi-norm-1}
\|f\|_{\omega,K,l}:=\sup_{k\in\NN,x\in K}\frac{\|f^{(k)}(x)\|_{L^k(\RR^r,\RR^s)}}{\exp(\frac{1}{l}\varphi^{*}_{\omega}(lk))}.
\end{equation}
For compact sets $K$ with smooth boundary
$$\mathcal{E}_{\omega,l}(K,\RR^s):=\{f\in\mathcal{E}(K,\RR^s): \|f\|_{\omega,K,l}<+\infty\}$$
is a Banach space and we consider the following topological vector spaces
\begin{equation}\label{repr1}
\mathcal{E}_{\{\omega\}}(U,\RR^s):=\underset{K\subset\subset U}{\varprojlim}\;\underset{l>0}{\varinjlim}\;\mathcal{E}_{\omega,l}(K,\RR^s)=\underset{K\subset\subset U}{\varprojlim}\;\mathcal{E}_{\{\omega\}}(K,\RR^s)
\end{equation}
and
\begin{equation}\label{repr2}
\mathcal{E}_{(\omega)}(U,\RR^s):=\underset{K\subset\subset U}{\varprojlim}\;\underset{l>0}{\varprojlim}\;\mathcal{E}_{\omega,l}(K,\RR^s)=\underset{K\subset\subset U}{\varprojlim}\;\mathcal{E}_{(\omega)}(K,\RR^s).
\end{equation}
For $\sigma,\tau\in\hyperlink{omset1}{\mathcal{W}}$ we get $\sigma\hyperlink{omprece}{\preceq}\tau\Leftrightarrow\mathcal{E}_{[\sigma]}\subseteq\mathcal{E}_{[\tau]}$ and $\tau\hyperlink{omtriangle}{\vartriangleleft}\sigma\Leftrightarrow\mathcal{E}_{\{\tau\}}\subseteq\mathcal{E}_{(\sigma)}$, see \cite[Corollary 5.17]{compositionpaper}.

We summarize some facts which are shown in \cite[Section 5]{compositionpaper}.
\begin{itemize}
\item[$(i)$] A central new idea was that to each $\omega\in\hyperlink{omset1}{\mathcal{W}}$ we can associate a \hyperlink{Msc}{$(\mathcal{M}_{\text{sc}})$} weight matrix $\Omega:=\{\Omega^l=(\Omega^l_j)_{j\in\NN}: l>0\}$ by

    \centerline{\fbox{$\Omega^l_j:=\exp\left(\frac{1}{l}\varphi^{*}_{\omega}(lj)\right)$.}}

\item[$(ii)$] $\mathcal{E}_{[\omega]}=\mathcal{E}_{[\Omega]}$ holds as locally convex vector spaces and $\Omega$ satisfies \hyperlink{R-mg}{$(\mathcal{M}_{\{\text{mg}\}})$}, \hyperlink{B-mg}{$(\mathcal{M}_{(\text{mg})})$} and \hyperlink{R-L}{$(\mathcal{M}_{\{\text{L}\}})$}, \hyperlink{B-L}{$(\mathcal{M}_{(\text{L})})$}.
\item[$(iii)$] Equivalent weight functions $\omega$ yield equivalent weight matrices w.r.t. both \hyperlink{Mbeurapprox}{$(\approx)$} and \hyperlink{Mroumapprox}{$\{\approx\}$}. Note that $(\mathcal{M}_{[\text{mg}]})$ is stable w.r.t. $[\approx]$, whereas $(\mathcal{M}_{[\text{L}]})$ not.
\item[$(iv)$] Defining classes of ultradifferentiable functions by weight matrices as in \eqref{generalroumieu} and in \eqref{generalbeurling} is a common generalization of defining them by using a (single) weight sequence $M$, i.e. a constant weight matrix, or by a weight function $\omega\in\hyperlink{omset1}{\mathcal{W}}$. But one is able to describe also other classes, e.g. the class defined by the Gevrey-matrix $\mathcal{G}:=\{(p!^{s+1})_{p\in\NN}: s>0\}$.
\end{itemize}

\subsection{Classes of ultra-differentiable functions defined by a weight matrix of associated functions}
Let $M\in\RR_{>0}^{\NN}$, the {\itshape associated function} $\omega_M: \RR_{\ge 0}\rightarrow\RR\cup\{+\infty\}$ is defined by
\begin{equation}\label{assofunc}
\omega_M(t):=\sup_{p\in\NN}\log\left(\frac{t^p M_0}{M_p}\right)\;\;\;\text{for}\;t>0,\hspace{30pt}\omega_M(0):=0.
\end{equation}

\begin{lemma}\label{assofuncproper}
If $M\in\hyperlink{LCset}{\mathcal{LC}}$, then $\omega_M$ belongs to \hyperlink{omset0}{$\mathcal{W}_0$}.

Moreover $\liminf_{p\rightarrow\infty}(m_p)^{1/p}>0$ implies \hyperlink{om2}{$(\omega_2)$}, $\lim_{p\rightarrow\infty}(m_p)^{1/p}=+\infty$ implies \hyperlink{om5}{$(\omega_5)$} for $\omega_M$.
\end{lemma}
We refer to \cite[Definition 3.1]{Komatsu73} and \cite[Lemma 12 $(iv)\Rightarrow(v)$]{BonetMeiseMelikhov07}. That $\lim(m_p)^{1/p}=+\infty$ implies \hyperlink{om5}{$(\omega_5)$} for $\omega_M$ follows analogously as $\liminf(m_p)^{1/p}>0$ implies \hyperlink{om2}{$(\omega_2)$} for $\omega_M$ as shown in \cite[Lemma 12 $(iv)\Rightarrow(v)$]{BonetMeiseMelikhov07}. Note that by {\itshape Stirling's formula} $\liminf(m_p)^{1/p}>0$ is precisely $(M0)$ in \cite{BonetMeiseMelikhov07}.

\begin{remark}\label{importantremark}
Let $\omega\in\hyperlink{omset0}{\mathcal{W}_0}$ be given, then
\begin{itemize}
\item[$(1)$] $\Omega^l\in\hyperlink{LCset}{\mathcal{LC}}$ for each $l>0$ by \cite[5.5]{compositionpaper},
\item[$(2)$] $\omega\hyperlink{sim}{\sim}\omega_{\Omega^l}$ for each $l>0$ by \cite[Lemma 5.7]{compositionpaper},
\item[$(3)$] $\omega$ satisfies
 \begin{itemize}
 \item[$(a)$] \hyperlink{omnq}{$(\omega_{\on{nq}})$} if and only if some/each $\Omega^l$ satisfies \hyperlink{mnq}{$(\on{nq})$},
 \item[$(b)$] \hyperlink{om6}{$(\omega_6)$} if and only if some/each $\Omega^l$ satisfies \hyperlink{mg}{$(\on{mg})$} if and only if $\Omega^l\hyperlink{approx}{\approx}\Omega^n$ for each $l,n>0$,
 \end{itemize}
 by \cite[Corollary 5.8, Theorem 5.14]{compositionpaper}.
\end{itemize}
\end{remark}
Let $\mathcal{M}=\{M^x: x\in\Lambda\}$ be \hyperlink{Msc}{$(\mathcal{M}_{\on{sc}})$}, then we introduce the new weight matrix $\omega_{\mathcal{M}}:=\{\omega_{M^x}: x\in\Lambda\}$. Let $U\subseteq\RR^r$ be non-empty open and put
$$\mathcal{E}_{\{\omega_{\mathcal{M}}\}}(U,\RR^s):=\{f\in\mathcal{E}(U,\RR^s):\;\forall\;K\subset\subset U\;\exists\;x\in\Lambda\;\exists\;l>0:\;\|f\|_{\omega_{M^x},K,l}<+\infty\}$$
and
$$\mathcal{E}_{(\omega_{\mathcal{M}})}(U,\RR^s):=\{f\in\mathcal{E}(U,\RR^s):\;\forall\;K\subset\subset U\;\forall\;x\in\Lambda\;\forall\;l>0:\;\|f\|_{\omega_{M^x},K,l}<+\infty\}.$$
Thus we obtain the topological vector spaces representations
\begin{equation}\label{assofuncroummatrix}
\mathcal{E}_{\{\omega_{\mathcal{M}}\}}(U,\RR^s):=\underset{K\subset\subset U}{\varprojlim}\;\underset{x\in\Lambda,l>0}{\varinjlim}\mathcal{E}_{\omega_{M^x},l}(K,\RR^s)
\end{equation}
and
\begin{equation}\label{assofuncbeurmatrix}
\mathcal{E}_{(\omega_{\mathcal{M}})}(U,\RR^s):=\underset{K\subset\subset U}{\varprojlim}\;\underset{x\in\Lambda,l>0}{\varprojlim}\mathcal{E}_{\omega_{M^x},l}(K,\RR^s)
\end{equation}

\section{Stability of constructing multi-index weight matrices}\label{section3}
\subsection{Introduction}\label{subsection30}
Let $\mathcal{M}:=\{M^x: x\in\Lambda\}$ be \hyperlink{Msc}{$(\mathcal{M}_{\text{sc}})$}. By Lemma \ref{assofuncproper} we get $\omega_{M^x}\in\hyperlink{omset0}{\mathcal{W}_0}$ for each $x\in\Lambda$. On the other hand by \cite[5.5]{compositionpaper} to each $\omega\in\hyperlink{omset0}{\mathcal{W}_0}$ we can associate a \hyperlink{Msc}{$(\mathcal{M}_{\text{sc}})$} weight matrix $\Omega:=\{(\Omega^l_j)_{j\in\NN}: l>0\}$ by putting $\Omega^l_j:=\exp\left(\frac{1}{l}\varphi^{*}_{\omega}(lj)\right)$.

So one can consider the construction
\begin{equation}\label{naturalitydefequation}
M^x\mapsto\omega_{M^x}\mapsto M^{x;l_1}\mapsto\omega_{M^{x;l_1}}\mapsto\ M^{x;l_1,l_2}\mapsto\dots,
\end{equation}
where for $x\in\Lambda$, $l_j\in\RR_{>0}$, $j\in\NN_{>0}$, and $i\in\NN$ we put
$$M^{x;l_1,\dots,l_{j+1}}_i:=\exp\left(\frac{1}{l_{j+1}}\varphi^{*}_{\omega^{M^{x;l_1,\dots,l_j}}}(l_{j+1}i)\right),\;\;\; M^{x;l_1}_i:=\exp\left(\frac{1}{l_1}\varphi^{*}_{\omega^{M^x}}(l_1i)\right)$$
respectively
$$\omega_{M^{x;l_1,\dots,l_j}}(t):=\sup_{p\in\NN}\log\left(\frac{t^p}{M_p^{x;l_1,\dots,l_j}}\right)\;\text{for}\;t>0,\;\hspace{20pt}\omega_{M^{x;l_1,\dots,l_j}}(0):=0.$$
On the one hand we obtain a sequence of matrices of weight functions. \cite[Lemma 5.7]{compositionpaper} implies
\begin{equation}\label{superequ2}
\forall\;x\in\Lambda\;\forall\;j\in\NN_{>0}\;\forall\;l_1,\dots,l_j>0:\;\omega_{M^{x;l_1,\dots,l_{j+1}}}\hyperlink{sim}{\sim}\omega_{M^{x;l_1,\dots,l_j}}\hyperlink{sim}{\sim}\dots\hyperlink{sim}{\sim}\omega_{M^x},
\end{equation}
hence this construction is always stable. So for each non-empty open $U\subseteq\RR^r$ we get
\begin{equation}\label{weightfunctionmatrix1}
\mathcal{E}_{\{\omega_{\mathcal{M}}\}}(U,\RR^s)=\underset{K\subset\subset U}{\varprojlim}\;\underset{x\in\Lambda,l,h>0}{\varinjlim}\mathcal{E}_{\omega_{M^{x;l}},h}(K,\RR^s)
\end{equation}
and
\begin{equation}\label{weightfunctionmatrix2}
\mathcal{E}_{(\omega_{\mathcal{M}})}(U,\RR^s)=\underset{K\subset\subset U}{\varprojlim}\;\underset{x\in\Lambda,l,h>0}{\varprojlim}\mathcal{E}_{\omega_{M^{x;l}},h}(K,\RR^s).
\end{equation}
On the other hand we get a sequence of matrices of weight sequences. In Theorem \ref{stablefirstcentralresult} we are going to characterize the stability of this construction and we will see that only in the first step of \eqref{naturalitydefequation} there can occur a non-stable effect (see Corollary \ref{1consequ}).

Finally the aim of this Section is to prove the following result:

\begin{theorem}\label{4}
Let $\mathcal{M}:=\{M^x\in\RR_{>0}^{\NN}: x\in\Lambda\}$ be \hyperlink{Msc}{$(\mathcal{M}_{\on{sc}})$}, let $r,s\in\NN_{>0}$ and $U$ be a non-empty open set in $\RR^r$. If $\mathcal{M}$ has $(\mathcal{M}_{[\on{L}]})$ and $(\mathcal{M}_{[\on{mg}]})$, then we get as locally convex vector spaces
$$\mathcal{E}_{[\mathcal{M}]}(U,\RR^s)=\mathcal{E}_{[\omega_{\mathcal{M}}]}(U,\RR^s).$$
\end{theorem}

\subsection{Stability of constructing multi-index matrices consisting of weight sequences}\label{subsection31}
In this section we show the following result which is the first step to prove Theorem \ref{4}.

\begin{theorem}\label{stablefirstcentralresult}
Let $\mathcal{M}=\{M^x: x\in\Lambda\}$ be \hyperlink{Msc}{$(\mathcal{M}_{\on{sc}})$}. Then $\mathcal{M}[\approx]\{M^{x;l}: x\in\Lambda, l>0\}$ if and only if
\begin{itemize}
\item[$(1)$] in the Roumieu-case \hyperlink{R-mg}{$(\mathcal{M}_{\{\on{mg}\}})$} holds,

\item[$(2)$] in the Beurling-case \hyperlink{B-mg}{$(\mathcal{M}_{(\on{mg})})$} holds, provided $\Lambda=\RR_{>0}$.
\end{itemize}
\end{theorem}

First we prove
\begin{lemma}
For each $x\in\Lambda$, $l\in\NN_{>0}$ and $j\in\NN$ we get
\begin{equation}\label{superequ1}
M_j^{x;l}=(M^x_{jl})^{1/l}.
\end{equation}
\end{lemma}

\demo{Proof}
We use \cite[Proposition 3.2]{Komatsu73} and get
\begin{align*}
M^{x;l}_j&:=\exp\left(\frac{1}{l}\varphi^{*}_{\omega_{M^x}}(lj)\right)=\exp\left(\frac{1}{l}\sup_{y\ge 0}\{y(lj)-\varphi_{\omega_{M^x}}(y)\}\right)=\exp\left(\sup_{y\ge 0}\left\{(y j)-\frac{1}{l}\varphi_{\omega_{M^x}}(y)\right\}\right)
\\&
=\sup_{y\ge 0}\frac{\exp(yj)}{\exp\left(\frac{1}{l}\varphi_{\omega_{M^x}}(y)\right)}=\sup_{s\ge 1}\frac{s^j}{\exp\left(\frac{1}{l}\omega_{M^x}(s)\right)}=\left(\sup_{s\ge 0}\frac{s^{jl}}{\exp(\omega_{M^x}(s))}\right)^{1/l}=(M^x_{jl})^{1/l}.
\end{align*}
All steps except the last one hold also for $l>0$ instead of $l\in\NN_{>0}$.
\qed\enddemo

The next result generalizes \cite[Proposition 3.6]{Komatsu73}.

\begin{proposition}\label{pseudomg}
Let $\mathcal{M}$ be \hyperlink{Msc}{$(\mathcal{M}_{\on{sc}})$}, then
\begin{equation}\label{pseudomg1}
\hyperlink{R-mg}{(\mathcal{M}_{\{\on{mg}\}})}\Longleftrightarrow\;\forall\;x\in\Lambda\;\exists\;H\ge 1\;\exists\;y\in\Lambda\;\forall\;t\ge 0: 2\omega_{M^y}(t)\le\omega_{M^x}(Ht)+H,
\end{equation}
\begin{equation}\label{pseudomg2}
\hyperlink{B-mg}{(\mathcal{M}_{(\on{mg})})}\Longleftrightarrow\;\forall\;x\in\Lambda\;\exists\;H\ge 1\;\exists\;y\in\Lambda\;\forall\;t\ge 0: 2\omega_{M^x}(t)\le\omega_{M^y}(Ht)+H.
\end{equation}
Even if $\omega_{M^x}\hyperlink{sim}{\sim}\omega_{M^y}$ for all $x,y\in\Lambda$, \eqref{pseudomg1} or \eqref{pseudomg2} does not imply necessarily \hyperlink{om6}{$(\omega_6)$} for each $\omega_{M^x}$.
\end{proposition}

\demo{Proof}
We follow \cite[Lemma 3.5, Proposition 3.6]{Komatsu73} and consider the Roumieu case. \hyperlink{R-mg}{$(\mathcal{M}_{\{\text{mg}\}})$} is equivalent to
$$\forall\;x\in\Lambda\;\exists\;H\ge 1\;\exists\;y\in\Lambda\;\forall\;p\in\NN:\;M^x_p\le H^p\min_{0\le q\le p} M^y_qM^y_{p-q}=:H^pN^y_p.$$
By \cite[Lemma 3.5]{Komatsu73} we have $\omega_{N^y}=2\omega_{M^y}$ and proceed as in \cite[Proposition 3.6]{Komatsu73} to get
\begin{align*}
2\omega_{M^y}(t)&=\sup_{p\in\NN}\log\left(\frac{t^p}{N^y_p}\right)=\sup_{p\in\NN}\log\left(\frac{t^p}{\min_{0\le q\le p} M^y_qM^y_{p-q}}\right)\le\sup_{p\in\NN}\log\left(\frac{t^pH^p}{M^x_p}\right)=\omega_{M^x}(Ht).
\end{align*}
Conversely, again as in \cite[Proposition 3.6]{Komatsu73}
\begin{align*}
N^y_p&=\sup_{t\ge 0}\frac{t^p}{\exp(\omega_{N^y}(t))}=\sup_{t\ge 0}\frac{t^p}{\exp(2\omega_{M^y}(t))}\ge\sup_{t\ge 0}\frac{t^p}{\exp(\omega_{M^x}(Ht)+H)}=\frac{1}{H^p\exp(H)}M^x_p.
\end{align*}
\qed\enddemo

Now we are able to prove the first part of Theorem \ref{stablefirstcentralresult}.

\begin{theorem}\label{1}
Let $\mathcal{M}:=\{M^x\in\RR_{>0}^{\NN}: x\in\Lambda\}$ be \hyperlink{Msc}{$(\mathcal{M}_{\on{sc}})$}, $r,s\in\NN_{>0}$. If \hyperlink{R-mg}{$(\mathcal{M}_{\{\on{mg}\}})$} holds then for each non-empty open set $U\subseteq\RR^r$ we get as locally convex vector spaces
$$\mathcal{E}_{\{\mathcal{M}\}}(U,\RR^s)=\underset{K\subset\subset U}{\varprojlim}\;\underset{x\in\Lambda,l,h>0}{\varinjlim}\mathcal{E}_{M^{x;l},h}(K,\RR^s).$$
If \hyperlink{B-mg}{$(\mathcal{M}_{(\on{mg})})$} holds then we get as locally convex vector spaces
$$\mathcal{E}_{(\mathcal{M})}(U,\RR^s)=\underset{K\subset\subset U}{\varprojlim}\;\underset{x\in\Lambda,l,h>0}{\varprojlim}\mathcal{E}_{M^{x;l},h}(K,\RR^s).$$
\end{theorem}

\demo{Proof}
{\itshape Roumieu case.} By \eqref{superequ1} implication $(\subseteq)$ holds in any case since $M^{x;1}=M^x\le M^{y}$ for $x\le y$. We show $(\supseteq)$ and by \eqref{superequ1} it suffices to prove
\begin{equation}\label{1-equ1}
\forall\;x\in\Lambda\;\forall\;l\in\NN_{>0}\;\exists\;y\in\Lambda\;\exists\;C\ge 1\;\forall\;j\in\NN:\;\;(M^x_{jl})^{1/l}\le C^jM^y_j\Leftrightarrow M^x_{jl}\le C^{jl}(M^y_j)^l,
\end{equation}
which implies $\mathcal{E}_{M^{x;l},h}(K,\RR^s)\subseteq\mathcal{E}_{M^y,Ch}(K,\RR^s)$. Now for each $x\in\Lambda$ there exists $D\ge 1$ and $y\in\Lambda$ such that $M^x_{2j}\le D^{2j}(M_j^{y})^2$ for all $j\in\NN$ by \hyperlink{R-mg}{$(\mathcal{M}_{\{\text{mg}\}})$} and so \eqref{1-equ1} follows by iterating this estimate $l$-times.\vspace{6pt}

{\itshape Beurling case.} $(\supseteq)$ is valid in any case since $M^{x;1}=M^x$ for each $x\in\Lambda$. Let us prove $(\subseteq)$, more precisely we show
\begin{equation}\label{1-equ2}
\forall\;x\in\Lambda\;\forall\;l>0\;\exists\;y\in\Lambda\;\exists\;C\ge 1\;\forall\;j\in\NN:\;\;M^y_j\le C^jM^{x;l}_j,
\end{equation}
which implies $\mathcal{E}_{M^y,h}(K,\RR^s)\subseteq\mathcal{E}_{M^{x;l},Ch}(K,\RR^s)$. Iterating \eqref{pseudomg2} gives
\begin{equation}\label{iteratedmgequ}
\forall\;x\in\Lambda\;\forall\;k\in\NN_{>0}\;\exists\;y\in\Lambda\;\exists\;H\ge 1\;\forall\;t\ge 0:\;2^k\omega_{M^x}(t)\le\omega_{M^y}(H^kt)+(2^k-1)H.
\end{equation}
Let $l\in\NN_{>0}$ be given (large) and $k\in\NN_{>0}$ be chosen minimal with $l\le 2^k$. For all $x\in\Lambda$ and $j\in\NN$ we have as in the proof of \eqref{superequ1}
\begin{align*}
M^{x;1/l}_j&=\sup_{t\ge 0}\frac{t^j}{\exp(l\omega_{M^x}(t))}
\ge\sup_{t\ge 0}\frac{t^j}{\exp(\omega_{M^{y}}(H^kt)+(2^k-1)H)}=\frac{1}{\exp((2^k-1)H)}\left(\frac{1}{H^k}\right)^jM^{y}_j.
\end{align*}
Consequently for arbitrary $x\in\Lambda$ and $l\in\NN_{>0}$ we find $y\in\Lambda$ such that $M^{y}\hyperlink{mpreceq}{\preceq}M^{x;1/l}$.
\qed\enddemo

An immediate consequence of Theorem \ref{1} is

\begin{corollary}\label{1consequ}
Let $\mathcal{M}$ be \hyperlink{Msc}{$(\mathcal{M}_{\on{sc}})$}, then after the first step in \eqref{naturalitydefequation} the construction yields always equivalent weight matrices of weight sequences w.r.t. to both \hyperlink{Mroumapprox}{$\{\approx\}$} and \hyperlink{Mbeurapprox}{$(\approx)$}.
\end{corollary}

\demo{Proof}
Let $x\in\Lambda$ be arbitrary but fixed. By Lemma \ref{assofuncproper} we have $\omega_{M^x}\in\hyperlink{omset0}{\mathcal{W}_0}$ and so \cite[5.5]{compositionpaper} implies that each matrix $\mathcal{M}^x:=\{M^{x;l}: l>0\}$, $x\in\Lambda$, satisfies both \hyperlink{R-mg}{$(\mathcal{M}_{\{\on{mg}\}})$} and \hyperlink{B-mg}{$(\mathcal{M}_{(\on{mg})})$}.
\qed\enddemo

Now we prove the converse implication for Theorem \ref{stablefirstcentralresult}. Here, the assumption $\Lambda=\RR_{>0}$ for the Beurling case is necessary.

\begin{proposition}\label{1converse}
Let $\mathcal{M}:=\{M^x\in\RR_{>0}^{\NN}: x\in\Lambda\}$ be \hyperlink{Msc}{$(\mathcal{M}_{\on{sc}})$}.
\begin{itemize}
\item[$(i)$] The equality
$$\mathcal{E}_{\{\mathcal{M}\}}(\RR,\RR)=\underset{K\subset\subset\RR}{\varprojlim}\;\underset{x\in\Lambda,l,h>0}{\varinjlim}\mathcal{E}_{M^{x;l},h}(K,\RR)$$
implies \hyperlink{R-mg}{$(\mathcal{M}_{\{\on{mg}\}})$} for $\mathcal{M}$.
\item[$(ii)$] Assume that $\Lambda=\RR_{>0}$, then
$$\mathcal{E}_{(\mathcal{M})}(\RR,\RR)=\underset{K\subset\subset\RR}{\varprojlim}\;\underset{x\in\Lambda,l,h>0}{\varprojlim}\mathcal{E}_{M^{x;l},h}(K,\RR)$$
implies \hyperlink{B-mg}{$(\mathcal{M}_{(\on{mg})})$} for $\mathcal{M}$.
\end{itemize}
\end{proposition}

\demo{Proof}
We generalize the technique in the proof of \cite[Lemma 5.9 $(5.11)$]{compositionpaper}.

{\itshape Roumieu case.} For each $x\in\Lambda$ and $l>0$ there exists a characteristic function $\theta_{x,l}\in\mathcal{E}_{\{M^{x;l}\}}^{\text{global}}(\RR,\RR)$, see \eqref{characteristicfunction}. So the inclusion $(\supseteq)$ implies
$$\forall\;x\in\Lambda\;\forall\;l>0\;\exists\;y\in\Lambda:\;M^{x;l}\hyperlink{mpreceq}{\preceq}M^{y}=M^{y;1},$$
equivalently
\begin{equation}\label{1converseequ1}
\forall\;x\in\Lambda\;\forall\;l>0\;\exists\;y\in\Lambda\;\exists\;C\ge 1\;\forall\;j\in\NN:\;\;\frac{1}{l}\varphi^{*}_{\omega_{M^{x}}}(lj)\le j\log(C)+\varphi^{*}_{\omega_{M^{y}}}(j).
\end{equation}
Consider \eqref{1converseequ1} for all $t\ge 0$ instead of all $j\in\NN$. Then
\begin{align*}
\left(\frac{1}{l}\varphi^{*}_{\omega_{M^{x}}}(l\cdot)\right)^{*}(s)&=\sup_{t\ge 0}\left\{st-\frac{1}{l}\varphi^{*}_{\omega_{M^{x}}}(lt)\right\}=\frac{1}{l}\sup_{t'\ge 0}\left\{st'-\varphi^{*}_{\omega_{M^{x}}}(t')\right\}
\\&
=\frac{1}{l}\varphi^{**}_{\omega_{M^{x}}}(s)=\frac{1}{l}\varphi_{\omega_{M^{x}}}(s)=\frac{1}{l}\omega_{M^{x}}(\exp(s)),
\end{align*}
which holds since $\omega_{M^x}\in\hyperlink{omset0}{\mathcal{W}_0}$ and so $\varphi^{**}_{\omega_{M^{x}}}(s)=\varphi_{\omega_{M^{x}}}(s)$. The right hand side gives
\begin{align*}
(\cdot D+\varphi^{*}_{\omega_{M^{y}}}(\cdot))^{*}(s)&=\sup_{t\ge 0}\{(s-D)t-\varphi^{*}_{\omega_{M^{y}}}(t)\}=\varphi^{**}_{\omega_{M^{y}}}(s-D)=\varphi_{\omega_{M^{y}}}(s-D)=\omega_{M^{y}}\left(\frac{\exp(s)}{C}\right).
\end{align*}
Then we use \cite[Lemma 5.7]{compositionpaper} (since $\omega_{M^x}\in\hyperlink{omset0}{\mathcal{W}_0}$ we can replace $\omega$ by $\omega_{M^{y}}=\omega_{M^{y;1}}$ there) and get for $s\ge 0$ sufficiently large:
\begin{align*}
&\sup_{t\ge 0}\left\{st-\frac{1}{l}\varphi^{*}_{\omega_{M^{x}}}(lt)\right\}\ge\sup_{j\in\NN}\left\{sj-\frac{1}{l}\varphi^{*}_{\omega_{M^{x}}}(lj)\right\}\underbrace{\ge}_{\eqref{1converseequ1}}\sup_{j\in\NN}\{sj-jD-\varphi^{*}_{\omega_{M^{y}}}(j)\}
\\&
\ge\frac{1}{2}\sup_{t\ge 0}\{st-tD-\varphi^{*}_{\omega_{M^{y}}}(t)\}=\frac{1}{2}\varphi^{**}_{\omega_{M^{y}}}(s-D)=\frac{1}{2}\omega_{M^{y}}\left(\frac{\exp(s)}{C}\right).
\end{align*}
Thus for all $t$ sufficiently large $\frac{1}{l}\omega_{M^{x}}(t)\ge\frac{1}{2}\omega_{M^{y}}\left(\frac{t}{C}\right)$ holds. Put $l=4$ and by \eqref{pseudomg1} we have shown \hyperlink{R-mg}{$(\mathcal{M}_{\{\on{mg}\}})$}.\vspace{6pt}

{\itshape Beurling case.} We follow the second Section in \cite{Bruna}, see also \cite[Proposition 4.6 $(1)$]{compositionpaper}. By assumption $\bigcap_{x\in\Lambda}\mathcal{E}_{(M^x)}(\RR,\RR)\subseteq\bigcap_{x\in\Lambda,l>0}\mathcal{E}_{(M^{x;l})}(\RR,\RR)$ and both are Fréchet spaces. Using the closed graph theorem the inclusion is continuous. Hence for each compact set $K_1\subseteq\RR$, $x\in\Lambda$, $l>0$ and $h>0$, there exist $C,h_1>0$, $y\in\Lambda$ and a compact set $K_2\subseteq\RR$ such that for each $f\in\bigcap_{x\in\Lambda}\mathcal{E}_{(M^x)}(\RR,\RR)$ we obtain
$$\|f\|_{M^{x;l},K_1,h}=\sup_{t\in K_1, j\in\NN}\frac{|f^{(j)}(t)|}{h^j M^{x;l}_j}\le C\sup_{t\in K_2, j\in\NN}\frac{|f^{(j)}(t)|}{h_1^jM^{y}_j}=C\|f\|_{M^y,K_2,h_1}.$$
Let $K_1$ be a compact interval containing $0$, put $h=1$ and take $f_s(t):=\sin(st)+\cos(st)$ for $t\in\RR$ and $s\ge 0$. Note that $f_s\in\bigcap_{x\in\Lambda}\mathcal{E}_{(M^x)}^{\text{global}}(\RR,\RR)$ for any $s\ge 0$ since $\lim_{k\rightarrow\infty}(M^x_k)^{1/k}=+\infty$ for each $x\in\Lambda$. Then
\begin{align*}
\sup_{j\in\NN}\frac{s^j}{M^{x;l}_j}&=\sup_{j\in\NN}\frac{|f_s^{(j)}(0)|}{M^{x;l}_j}\le\sup_{t\in K_1, j\in\NN}\frac{|f_s^{(j)}(t)|}{M^{x;l}_j}\le C\sup_{t\in K_2, j\in\NN}\frac{|f_s^{(j)}(t)|}{h_1^j M^{y}_j}\le C\sup_{j\in\NN}\frac{2 s^j}{h_1^jM^{y}_j},
\end{align*}
which implies $\exp(\omega_{M^{x;l}}(s))\le 2C\exp\left(\omega_{M^{y}}\left(\frac{s}{h_1}\right)\right)$. Using \cite[Proposition 3.2]{Komatsu73} we get for all $j\in\NN$
$$M^{x;l}_j=\sup_{t\ge 0}\frac{t^j}{\exp(\omega_{M^{x;l}}(t))}\ge\sup_{t\ge 0}\frac{t^j}{2C\exp\left(\omega_{M^{y}}\left(\frac{t}{h_1}\right)\right)}=\frac{h_1^j}{2C}M^{y}_j,$$
hence $M^{y}\hyperlink{mpreceq}{\preceq}M^{x;l}$. We summarize:
\begin{equation}\label{1converseequ2}
\forall\;x\in\Lambda\;\forall\;l>0\;\exists\;y\in\Lambda\;\exists\;D\ge 1\;\forall\;j\in\NN:\;\;\varphi^{*}_{\omega_{M^{y}}}(j)\le j\log(D)+\frac{1}{l}\varphi^{*}_{\omega_{M^{x}}}(lj).
\end{equation}
Now use the proof of the Roumieu case to get $\omega_{M^{y}}(t)\ge\frac{1}{2l}\omega_{M^{x}}\left(\frac{t}{D}\right)$ for $t$ sufficiently large. The choice $l=\frac{1}{4}$ and \eqref{pseudomg2} imply \hyperlink{B-mg}{$(\mathcal{M}_{(\on{mg})})$}.
\qed\enddemo

\subsection{Classes $\mathcal{E}_{[\omega_{\mathcal{M}}]}$ defined by a weight matrix of associated functions}\label{subsection32}
The goal of this section is to prove

\begin{theorem}\label{prä4}
Let $\mathcal{M}:=\{M^x\in\RR_{>0}^{\NN}: x\in\Lambda\}$ be \hyperlink{Msc}{$(\mathcal{M}_{\on{sc}})$}, let $r,s\in\NN_{>0}$ and $U$ be a non-empty open set in $\RR^r$.
\begin{itemize}
\item[$(i)$] \hyperlink{R-L}{$(\mathcal{M}_{\{\on{L}\}})$} for $\mathcal{M}$ implies
$$\mathcal{E}_{\{\omega_{\mathcal{M}}\}}(U,\RR^s)=\underset{K\subset\subset U}{\varprojlim}\underset{x\in\Lambda,l,h>0}{\varinjlim}\mathcal{E}_{M^{x;l},h}(K,\RR^s),$$
\item[$(ii)$] \hyperlink{B-L}{$(\mathcal{M}_{(\on{L})})$} for $\mathcal{M}$ implies
$$\mathcal{E}_{(\omega_{\mathcal{M}})}(U,\RR^s)=\underset{K\subset\subset U}{\varprojlim}\;\underset{x\in\Lambda,l,h>0}{\varprojlim}\mathcal{E}_{M^{x;l},h}(K,\RR^s)$$
\end{itemize}
as locally convex vector spaces.
\end{theorem}
The main Theorem \ref{4} follows then by combining Theorem \ref{1} and Theorem \ref{prä4}.

We start with the following result:

\begin{proposition}\label{2}
Let $\mathcal{M}:=\{M^x\in\RR_{>0}^{\NN}: x\in\Lambda\}$ be \hyperlink{Msc}{$(\mathcal{M}_{\on{sc}})$}.

\begin{itemize}
\item[$(i)$] \hyperlink{R-L}{$(\mathcal{M}_{\{\on{L}\}})$} implies
\begin{equation}\label{R-L-consequ}
\forall\;x\in\Lambda\;\exists\;y\in\Lambda:\;\;\;\omega_{M^y}(2t)=O(\omega_{M^x}(t))\;\;\text{as}\;t\rightarrow\infty.
\end{equation}
\item[$(ii)$] \hyperlink{B-L}{$(\mathcal{M}_{(\on{L})})$} implies
\begin{equation}\label{B-L-consequ}
\forall\;x\in\Lambda\;\exists\;y\in\Lambda:\;\;\;\omega_{M^x}(2t)=O(\omega_{M^y}(t))\;\;\text{as}\;t\rightarrow\infty.
\end{equation}
\end{itemize}
If all associated functions are equivalent w.r.t. $\hyperlink{sim}{\sim}$, then each/some $\omega_{M^x}$ satisfies \hyperlink{om1}{$(\omega_1)$}.
\end{proposition}

\demo{Proof}
By \hyperlink{R-L}{$(\mathcal{M}_{\{\text{L}\}})$} for each $x\in\Lambda$ and each $h>0$ there exists $y\in\Lambda$ and $D>0$ such that $M^x_kh^k\le D M^y_k$ holds for all $k\in\NN$. Multiplying with $t^k$ for arbitrary $t>0$ we get $\frac{(ht)^k}{M^y_k}\le D\frac{t^k}{M^x_k}$ and finally $\log\left(\frac{(ht)^k}{M^y_k}\right)\le\log\left(\frac{t^k}{M^x_k}\right)+D_1$,
which holds for all $k\in\NN$. So by definition $\omega_{M^y}(ht)\le\omega_{M^x}(t)+D_1$ holds and it is enough to take $h=2$.

The Beurling case is completely analogous, use \hyperlink{B-L}{$(\mathcal{M}_{(\text{L})})$} instead of \hyperlink{R-L}{$(\mathcal{M}_{\{\text{L}\}})$}.
\qed\enddemo

The next result generalizes \cite[Lemma 5.9 $(5.10)$]{compositionpaper}.

\begin{proposition}\label{3}
Let $\{\sigma_x\in\hyperlink{omset0}{\mathcal{W}_0}: x\in\Lambda\}$ be given and assume the Roumieu type condition (see Proposition \ref{2} above):
$$\forall\;x\in\Lambda\;\exists\;y\in\Lambda:\;\sigma_y(2t)=O(\sigma_x(t))\;\;\text{as}\;t\rightarrow\infty.$$
Then
\begin{align*}
&\forall\;x\in\Lambda\;\forall\;s\in\NN\;\exists\;y\in\Lambda\;\exists\;L\ge 1\;\forall\;a>0\;\forall\;j\in\NN:
\\&
\exp\left(\frac{1}{a}\varphi^{*}_{\sigma_x}(aj)\right)\exp(s)^j\le\exp\left(\frac{\sum_{i=1}^s L^i}{L^s a}\right)\exp\left(\frac{1}{L^sa}\varphi^{*}_{\sigma_y}(L^saj)\right).
\end{align*}
If the Beurling type condition
$$\forall\;x\in\Lambda\;\exists\;y\in\Lambda:\;\sigma_x(2t)=O(\sigma_y(t))\;\;\text{as}\;t\rightarrow\infty$$
holds, then
\begin{align*}
&\forall\;x\in\Lambda\;\forall\;s\in\NN\;\exists\;y\in\Lambda\;\exists\;L\ge 1\;\forall\;a>0\;\forall\;j\in\NN:
\\&
\exp\left(\frac{1}{a}\varphi^{*}_{\sigma_y}(aj)\right)\exp(s)^j\le\exp\left(\frac{\sum_{i=1}^s L^i}{L^s a}\right)\exp\left(\frac{1}{L^sa}\varphi^{*}_{\sigma_x}(L^saj)\right).
\end{align*}
If each $\omega_{M^x}$ has \hyperlink{om1}{$(\omega_1)$}, then the Roumieu and the Beurling case is satisfied with $x=y$.
\end{proposition}

\demo{Proof}
We consider the Roumieu case. For all $x\in\Lambda$ there exist $y\in\Lambda$ and $L\ge 1$ with $\sigma_y(4t)\le L\sigma_x(t)+L$ for all $t\ge 0$, hence $\varphi_{\sigma_y}(t+1)=\sigma_y(\exp(t+1))\le L\sigma_x(\exp(t))+L$. First we have
\begin{align*}
\varphi^{*}_{\sigma_y}(Ls)&=L\sup\left\{st-\frac{1}{L}\varphi_{\sigma_y}(t):t\ge 0\right\}\ge L\sup\{s t-(1+\varphi_{\sigma_x}(t-1)):t\ge0\}
\\&
\ge L\sup\{s(t-1)+s-1-\varphi_{\sigma_x}(t-1):t\ge1\}=Ls-L+L\varphi^{*}_{\sigma_x}(s),
\end{align*}
and so
$$\forall\;x\in\Lambda\;\exists\;y\in\Lambda\;\exists\;L\ge 1\;\forall\;t\ge 0:\;\;L\varphi^{*}_{\sigma_x}(t)+Lt\le L+\varphi^{*}_{\sigma_y}(Lt).$$
Using induction on this inequality we get
$$\forall\;x\in\Lambda\;\forall\;s\in\NN\;\exists\;y\in\Lambda\;\exists\;L\ge 1\;\forall\;t\ge 0:\;\;L^s\varphi^{*}_{\sigma_x}(t)+s L^st\le\varphi^{*}_{\sigma_y}(L^s t)+\sum_{i=1}^s L^i.$$
Now put $t=aj$ for $j\in\NN$ and $a>0$, divide by $L^sa$ and finally apply $\exp$.
\qed\enddemo

Propositions \ref{2} and \ref{3} imply

\begin{corollary}\label{naturalitycor}
Let $\mathcal{M}:=\{M^x\in\RR_{>0}^{\NN}: x\in\Lambda\}$ be \hyperlink{Msc}{$(\mathcal{M}_{\on{sc}})$}.
\begin{itemize}
\item[$(i)$] If $\mathcal{M}$ has \hyperlink{R-L}{$(\mathcal{M}_{\{\on{L}\}})$}, then
\begin{equation}\label{naturalitycor1}
\forall\;x\in\Lambda\;\forall\;h>0\;\exists\;y\in\Lambda\;\forall\;a>0\;\exists\;D>0\;\exists\;b>0\;\forall\;j\in\NN: M_j^{x;a}h^j\le D M_j^{y;b}.
\end{equation}
\item[$(ii)$] If $\mathcal{M}$ has \hyperlink{B-L}{$(\mathcal{M}_{(\on{L})})$}, then
\begin{equation}\label{naturalitycor2}
\forall\;x\in\Lambda\;\forall\;h>0\;\exists\;y\in\Lambda\;\forall\;b>0\;\exists\;D>0\;\exists\;a>0\;\forall\;j\in\NN: M_j^{y;a}h^j\le D M_j^{x;b}.
\end{equation}
\end{itemize}
\end{corollary}

Using \eqref{naturalitycor1} in the Roumieu and \eqref{naturalitycor2} in the Beurling case we get Theorem \ref{prä4} and are done.

We can also prove:

\begin{corollary}\label{naturalitycorchar}
Let $\mathcal{M}:=\{M^x\in\RR_{>0}^{\NN}: x\in\Lambda\}$ be \hyperlink{Msc}{$(\mathcal{M}_{\on{sc}})$}, then \eqref{R-L-consequ}$\Longleftrightarrow$\eqref{naturalitycor1} and \eqref{B-L-consequ}$\Longleftrightarrow$\eqref{naturalitycor2}.
\end{corollary}

\demo{Proof}
It remains to show $(\Longleftarrow)$. In \eqref{naturalitycor1} let $h=2$, $a=1$, multiply with $t^j$ for arbitrary $t>0$ and apply $\log$. Thus $\omega_{M^{y;b}}(2t)=O(\omega_{M^x}(t))$ holds as $t\rightarrow\infty$. Finally \cite[Lemma 5.7]{compositionpaper} implies $\omega_{M^{y;b}}\hyperlink{sim}{\sim}\omega_{M^y}$ for each $b>0$. The case for \eqref{naturalitycor2} is analogous.
\qed\enddemo

\subsection{Applications of Theorem \ref{4}}\label{subsection33}
If $\mathcal{M}=\Omega$ for some $\omega\in\hyperlink{omset1}{\mathcal{W}}$, then by Theorem \ref{4} and \cite[Theorem 5.14]{compositionpaper} we get $\mathcal{E}_{[\omega]}=\mathcal{E}_{[\Omega]}=\mathcal{E}_{[\omega_{\Omega}]}=\mathcal{E}_{[\omega_{\Omega^l}]}$ for each $l>0$. More generally we can prove

\begin{corollary}
Let $\mathcal{M}=\{M^x: x>0\}$ have \hyperlink{Msc}{$(\mathcal{M}_{\on{sc}})$}. Then the following are equivalent:

\begin{itemize}
\item[$(i)$] There exists $\omega\in\hyperlink{omset}{\mathcal{W}}$ with $\mathcal{E}_{[\mathcal{M}]}=\mathcal{E}_{[\omega]}$.
\item[$(ii)$] There exists a \hyperlink{Msc}{$(\mathcal{M}_{\on{sc}})$}-matrix $\mathcal{N}=\{N^x: x>0\}$ with $\mathcal{M}[\approx]\mathcal{N}$, such that $\omega_{N^x}\hyperlink{sim}{\sim}\omega_{N^y}$ for each $x,y>0$ and $\mathcal{N}$ has $(\mathcal{M}_{[\on{mg}]})$ and $(\mathcal{M}_{[\on{L}]})$.
\end{itemize}
\end{corollary}

\demo{Proof}
$(i)\Rightarrow(ii)$ We can take $\mathcal{N}=\Omega$, see \cite[Proposition 4.6, Lemma 5.7]{compositionpaper} and \cite[Theorem 5.14, Corollary 5.15]{compositionpaper}.

$(ii)\Rightarrow(i)$ Combining Theorem \ref{4} and \cite[Theorem 5.14]{compositionpaper} we get
\begin{equation}\label{charcequ1}
\forall\;x>0:\;\;\mathcal{E}_{[\mathcal{M}]}=\mathcal{E}_{[\mathcal{N}]}=\mathcal{E}_{[\omega_{\mathcal{N}}]}=\mathcal{E}_{[\omega_{N^x}]}=\mathcal{E}_{[\mathcal{N}^x]}
\end{equation}
with $\mathcal{N}^x:=\{N^{x;l}: l>0\}$. Note that $\omega_{M^x}\in\hyperlink{omset1}{\mathcal{W}}$ for each $x>0$, see Proposition \ref{2}. So we can take $\omega=\omega_{N^x}$ and $\Omega=\mathcal{N}^x$ for some arbitrary $x>0$, i.e. $\Omega^l=N^{x;l}$.

Finally by \cite[Proposition 4.6]{compositionpaper} we get $\mathcal{M}[\approx]\mathcal{N}^x$ and any $\sigma\in\hyperlink{omset1}{\mathcal{W}}$ with $\mathcal{E}_{[\sigma]}=\mathcal{E}_{[\mathcal{M}]}$ satisfies $\sigma\hyperlink{sim}{\sim}\omega_{N^x}$ by \cite[Corollary 5.17]{compositionpaper}.
\qed\enddemo

Let $\mathcal{M}=\{M^x: x\in\Lambda\}$ be \hyperlink{Msc}{$(\mathcal{M}_{\on{sc}})$} given, then in general we will not have $\omega_{M^x}\hyperlink{sim}{\sim}\omega_{M^y}$ for any $x,y\in\Lambda$. On the one hand by definition $\omega_{M^y}\le\omega_{M^x}$ whenever $x\le y$ and on the other hand \cite[1.8 III]{mandelbrojtbook} yields $\omega_{M^x}(t)=\sup_{p\in\NN}p\log(t)-\log(M^x_p)=p_{t,x}\log(t)-\log(M^x_{p_{t,x}})$, where $\mu^x_{p_{t,x}}\le t<\mu^x_{p_{t,x}+1}$. So if $\mathcal{M}$ satisfies
\begin{equation}\label{assofunccomparison}
\forall\;x,y>0\;x\le y\;\exists\;C\ge 1\;\exists\;t_0\ge 1\;\forall\;t\ge t_0\;\exists\;q\in\NN:\;\;\frac{(M^y_q)^C}{M^x_{p_{t,x}}}\le t^{qC-p_{t,x}},
\end{equation}
then all associated functions are equivalent w.r.t. \hyperlink{sim}{$\sim$}. Moreover we can prove:

\begin{lemma}\label{relationcomparison}
Let $M,N\in\hyperlink{LCset}{\mathcal{LC}}$.
\begin{itemize}
\item[$(1)$] If $\omega_M$ satisfies \hyperlink{om1}{$(\omega_1)$}, then $M\hyperlink{mpreceq}{\preceq} N\Longrightarrow\omega_M\hyperlink{ompreceq}{\preceq}\omega_N$.

\item[$(2)$] If $N$ satisfies \hyperlink{mg}{$(\on{mg})$}, then $\omega_N\hyperlink{ompreceq}{\preceq}\omega_M\Longrightarrow N\hyperlink{mpreceq}{\preceq}M$.
\end{itemize}
\end{lemma}

\demo{Proof}
$(1)$ For all $t>0$ we get
\begin{align*}
\omega_M(t)&=\sup_{p\in\NN}(p\log(t)-\log(M_p))\underbrace{\ge}_{M\hyperlink{mpreceq}{\preceq}N}\sup_{p\in\NN}(p\log(t)-\log(D^p N_p))=\omega_N\left(\frac{t}{D}\right)
\end{align*}
for a constant $D>0$ (large). Iterating \hyperlink{om1}{$(\omega_1)$} we have $\omega_M(2^n t)\le C^n\omega_M(t)+C$ for a constant $C\ge 1$ and all $t\ge 0$. Choose now $n\in\NN$ minimal such that $D\le 2^n$, hence $\omega_N(t)\le\omega_M(D t)\le\omega_M(2^n t)\le C^n\omega_M(t)+C$ for all $t\ge 0$ and so $\omega_N(t)=O(\omega_M(t))$ as $t\rightarrow\infty$.\vspace{6pt}

$(2)$ By \cite[Proposition 3.6]{Komatsu73} condition \hyperlink{mg}{$(\text{mg})$} for $N$ implies \hyperlink{om6}{$(\omega_6)$} for $\omega_N$. Using \cite[Proposition 3.2]{Komatsu73} we can estimate for all $p\in\NN$:
\begin{align*}
M_p&=\sup_{t>0}\frac{t^p}{\exp(\omega_M(t))}\underbrace{\ge}_{\omega_M\hyperlink{ompreceq}{\preceq}\omega_N}\sup_{t>0}\frac{t^p}{\exp(C_1\omega_N(t)+C_1)}\ge C_2\sup_{t>0}\frac{t^p}{\exp(\omega_N(H^nt)+(2^n-1)H)}
\\&
=C_3\left(\frac{1}{H^n}\right)^pN_p,
\end{align*}
where $n\in\NN$ is chosen minimal such that $C_1\le 2^n$ (iterating \hyperlink{om6}{$(\omega_6)$} as in \eqref{iteratedmgequ}). Thus $N\hyperlink{mpreceq}{\preceq}M$ follows.
\qed\enddemo

\subsection{Roumieu case versus Beurling case}\label{subsection34}
For $\mathcal{E}_{[\mathcal{M}]}$ and $\mathcal{E}_{[\omega_{\mathcal{M}}]}$ it is also important to know whether one can replace in their definitions the Roumieu classes $\mathcal{E}_{\{M^x\}}$, $\mathcal{E}_{\{\omega_{M^x}\}}$ by the Beurling classes $\mathcal{E}_{(M^x)}$, $\mathcal{E}_{(\omega_{M^x})}$. In the case $\mathcal{E}_{[\mathcal{M}]}$ this can be done assuming $(\mathcal{M}_{[\on{BR}]})$, see \cite[4.2 $(4.4)$]{compositionpaper}. If $\mathcal{M}=\Omega$ for some  $\omega\in\hyperlink{omset1}{\mathcal{W}}$, then \hyperlink{om7}{$(\omega_7)$} is sufficient to guarantee this property for the Roumieu case and the Beurling case, see \cite[Theorem 5.14 (4)]{compositionpaper}.

\begin{proposition}
Let $\mathcal{M}=\{M^x: x\in\Lambda\}$ be \hyperlink{Msc}{$(\mathcal{M}_{\on{sc}})$}.

\begin{itemize}
\item[$(i)$] If $\mathcal{M}$ has \hyperlink{R-BR}{$(\mathcal{M}_{\{\on{BR}\}})$} and each $M^x$ has \hyperlink{mg}{$(\on{mg})$}, then
\begin{equation}\label{triangle1}
\forall\;x\in\Lambda\;\exists\;y\in\Lambda:\;\;\omega_{M^x}\hyperlink{omtriangle}{\vartriangleleft}\omega_{M^y},
\end{equation}
which implies $\bigcup_{x\in\Lambda}\mathcal{E}_{\{\omega_{M^x}\}}=\bigcup_{x\in\Lambda}\mathcal{E}_{(\omega_{M^x})}$.

\item[$(ii)$] If $\mathcal{M}$ has \hyperlink{B-BR}{$(\mathcal{M}_{(\on{BR})})$} and each $M^x$ has \hyperlink{mg}{$(\on{mg})$}, then
\begin{equation}\label{triangle2}
\forall\;x\in\Lambda\;\exists\;y\in\Lambda:\;\;\omega_{M^y}\hyperlink{omtriangle}{\vartriangleleft}\omega_{M^x},
\end{equation}
which implies $\bigcap_{x\in\Lambda}\mathcal{E}_{\{\omega_{M^x}\}}=\bigcap_{x\in\Lambda}\mathcal{E}_{(\omega_{M^x})}$.

\item[$(iii)$] If each $\omega_{M^x}$ has \hyperlink{om1}{$(\omega_1)$} and \eqref{triangle1} holds, then $\mathcal{M}$ has \hyperlink{R-BR}{$(\mathcal{M}_{\{\on{BR}\}})$}.

\item[$(iv)$] If each $\omega_{M^x}$ has \hyperlink{om1}{$(\omega_1)$} and \eqref{triangle2} holds, then $\mathcal{M}$ has \hyperlink{B-BR}{$(\mathcal{M}_{(\on{BR})})$}.
\end{itemize}
\end{proposition}

\demo{Proof}
We consider the Roumieu case $(i)$ and $(iii)$, the Beurling case $(ii)$ and $(iv)$ is completely analogous.

$(i)$ \eqref{triangle1} means
$$\forall\;x\in\Lambda\;\exists\;y\in\Lambda\;\forall\;C>0\;\exists\;D>0\;\forall\;t\ge 0:\;\;\omega_{M^y}(t)\le C\omega_{M^x}(t)+D.$$
By assumption \hyperlink{R-BR}{$(\mathcal{M}_{\{\text{BR}\}})$} holds, i.e.
$$\forall\;x\in\Lambda\;\exists\;y\in\Lambda\;\forall\;h>0\;\exists\;C_h>0\;\forall\;j\in\NN:\;\;M^x_j\le C_h\ h^jM^y_j.$$
Multiplying with $t^j$ for arbitrary $t>0$ and $j\in\NN$ we get by definition $\log(C_h)+\omega_{M^x}(t)\ge\omega_{M^y}(t/h)$.

Now let $1>C>0$ be given, \hyperlink{mg}{$(\on{mg})$} for $M^y$ implies \hyperlink{om6}{$(\omega_6)$} for $\omega_{M^y}$. Iterating this condition (see \eqref{iteratedmgequ}) we take $k\in\NN$ minimal with $C^{-1}\le 2^k$ and choose $h:=\frac{1}{H^k}$. Then $C^{-1}\omega_{M^y}(t)\le\omega_{M^y}(H^kt)+(2^k-1)H=\omega_{M^y}(t/h)+H_1\le\omega_{M^x}(t)+H_2$.\vspace{6pt}

$(iii)$ Iterating \hyperlink{om1}{$(\omega_1)$} for $\omega_{M^x}$ gives $\omega_{M^x}(2^nt)\le L^n\omega_{M^x}(t)+\sum_{i=1}^n L^i$. So let $1>h>0$ be given and choose $n\in\NN_{>0}$ minimal with $h^{-1}\le 2^n$. Then $\omega_{M^x}(t/h)\le\omega_{M^x}(2^n t)\le L^n\omega_{M^x}(t)+\sum_{i=1}^n L^i$ and choose $C:=L^{-n}$ which depends only on $x\in\Lambda$ and given $h$. According to $x\in\Lambda$ and $C$ we use \eqref{triangle1} and \cite[Proposition 3.2]{Komatsu73} to obtain, for all $j\in\NN$:
\begin{align*}
M^y_j&=\sup_{t\ge 0}\frac{t^j}{\exp(\omega_{M^y}(t))}\ge\sup_{t\ge 0}\frac{t^j}{\exp(C\omega_{M^x}(t)+D)}\ge\frac{1}{D_1}\sup_{t\ge 0}\frac{t^j}{\exp(\omega_{M^x}(ht))}=\frac{1}{D_1 h^j}M^x_j.
\end{align*}
Note that the constant $D_1$ depends also only on $x$ and $h$.
\qed\enddemo

\section{Characterization of the non-quasianalyticity of $\mathcal{E}_{[\mathcal{M}]}$}\label{section4}
Let $\mathcal{M}$ be \hyperlink{Marb}{$(\mathcal{M})$}, then
$\mathcal{E}_{[\mathcal{M}]}$ is called {\itshape non-quasianalytic} if $\mathcal{E}_{[\mathcal{M}]}$ contains non-trivial functions with compact support.

The goal is to characterize this property in terms of the weight matrix $\mathcal{M}$ which gives answer to \cite[Remark 4.8]{compositionpaper}.

\begin{theorem}\label{Matrix-non-quasi-analyticity}
Let $\mathcal{M}=\{M^x: x\in\Lambda\}$ be \hyperlink{Marb}{$(\mathcal{M})$}.
\begin{itemize}
\item[$(i)$] $\mathcal{E}_{\{\mathcal{M}\}}$ is non-quasianalytic if and only if there exists $x_0\in\Lambda$ such that $\mathcal{E}_{[M^{x_0}]}$ is non-quasianalytic.

\item[$(ii)$] $\mathcal{E}_{(\mathcal{M})}$ is non-quasianalytic if and only if each $\mathcal{E}_{[M^{x}]}$ is non-quasianalytic, provided $\Lambda=\RR_{>0}$.
\end{itemize}
\end{theorem}

\begin{remark}
The theorem above still holds if we assume that each $M^x\in\RR^{\NN}_{>0}$ is arbitrary with $M^x_0=1$ and $M^x\le M^y$ whenever $x\le y$, i.e. the assumption that each $M^x$ is increasing is not necessary. This holds by the definitions of $\mathcal{E}_{[\mathcal{M}]}$ given in \ref{weighmatrixdefinition} and since we work in the proofs of Propositions \ref{nonquasiremarks} and \ref{Beurling-non-quasi-analyticity} below with the regularizations $M^{\on{lc}}$ and $M^I$ which will be defined in \ref{subsection41}. Note that $M\le N$ implies $M^{\on{lc}}\le N^{\on{lc}}$ and $M^I\le N^I$.
\end{remark}

\subsection{Non-quasianalyticity of $\mathcal{E}_{[M]}$}\label{subsection41}
Before we start proving Theorem \ref{Matrix-non-quasi-analyticity} we recall and summarize some facts for classical Denjoy-Carleman-classes $\mathcal{E}_{[M]}$. Let $M\in\RR_{>0}^{\NN}$ with $M_0=1$, then we denote by $M^{\on{lc}}=(M^{\on{lc}}_j)_j$ the log-convex minorant of $M$ which is given by
$$M^{\on{lc}}_j:=\sup_{t>0}\frac{t^j}{\exp(\omega_M(t))}$$
resp.
$$M^{\on{lc}}_j:=\inf\{M^{(l-j)/(l-k)}_k M^{(j-k)/(l-k)}: k\le j\le l, k\neq l\}, M^{\on{lc}}_0:=M_0=1$$
see \cite[Definition 3.1]{Komatsu73} and \cite{mandelbrojtbook} resp. \cite{hoermander}. Moreover we introduce
$$\hypertarget{bi}{M^I}:=(M^I_k)_k,\hspace{20pt}M^I_k:=\left(\inf\{(M_j)^{1/j}: j\ge k\}\right)^k\;\text{for}\;k\ge 1,\hspace{20pt} M^I_0:=1,$$
see also \cite{hoermander}. $((M^I_k)^{1/k})_k$ is the increasing minorant of $((M_k)^{1/k})_k$, $M^I=M$ if and only if $k\mapsto(M_k)^{1/k}$ is increasing. If $M$ is \hyperlink{lc}{$(\on{lc})$}, then $M=M^I$ and $(M^{\on{lc}})^I=M^{\on{lc}}$, so $M^{\on{lc}}\le M^I\le M$.

\begin{proposition}\label{nonquasiremarks}
Let $M\in\RR_{>0}^{\NN}$ with $M_0=1$. Then $\mathcal{E}_{[M]}$ is non-quasianalytic if and only if $M^{\on{lc}}$ has \hyperlink{mnq}{$(\on{nq})$} and if and only if $\sum_{p\ge 1}\frac{1}{(M^I_p)^{1/p}}<+\infty$. In this case $\mathcal{C}^{\omega}\subsetneq\mathcal{E}_{[M]}=\mathcal{E}_{[M^{I}]}=\mathcal{E}_{[M^{\on{lc}}]}$ holds.
\end{proposition}
{\itshape Remark:} The equivalence $\sum_{p=1}^{\infty}\frac{1}{(M^I_p)^{1/p}}<+\infty$ if and only if $M^{\on{lc}}$ has \hyperlink{mnq}{$(\on{nq})$} can be shown directly without using the non-quasianalyticity of $\mathcal{E}_{[M]}$, see the proof of \cite[Theorem 1.3.8]{hoermander}.

\demo{Proof}
By \cite[Theorem 1.3.8]{hoermander} and \cite[Theorem 4.2]{Komatsu73} we know that $\mathcal{E}_{[M]}$ is non-quasianalytic if and only if $\sum_{p=1}^{\infty}\frac{1}{(M^I_p)^{1/p}}<+\infty$ and if and only if \hyperlink{mnq}{$(\on{nq})$} holds for $M^{\on{lc}}$. More precisely the Roumieu-case follows directly by \cite[Theorem 1.3.8]{hoermander}. If $\mathcal{E}_{(M)}$ is non-quasianalytic, then $\mathcal{E}_{\{M\}}$ too and apply \cite[Theorem 1.3.8]{hoermander}. If $M^{\on{lc}}$ has \hyperlink{mnq}{$(\on{nq})$}, then by \cite[Theorem 4.2]{Komatsu73} the class $\mathcal{E}_{(M^{\on{lc}})}$ is non-quasianalytic, hence $\mathcal{E}_{(M)}$ too.\vspace{6pt}

{\itshape Claim.} If $\mathcal{E}_{[M]}$ is non-quasianalytic, then $\lim_{p\rightarrow\infty}(m_p)^{1/p}=+\infty\Leftrightarrow\lim_{k\rightarrow\infty}\frac{(M_p)^{1/p}}{p}=+\infty$. We put $a_p:=\frac{1}{(M^I_p)^{1/p}}$ in the well-known Lemma \ref{nonquasiremarkslemma} below and since $(M^I_p)^{1/p}\le(M_p)^{1/p}$ for all $p\in\NN_{>0}$ the claim follows.

This claim generalizes remark $(b.1)$ on page $387$ in \cite{surjectivity} since there \hyperlink{lc}{$(\on{lc})$} (which is assumed in $(b)$) for $M$ was necessary. Moreover it implies $\mathcal{C}^{\omega}\subsetneq\mathcal{E}_{[M]}$.

Finally by \cite[Theorem 2.15]{compositionpaper} and the claim we see that $\mathcal{E}_{[M^{\on{lc}}]}\subseteq\mathcal{E}_{[M^I]}\subseteq\mathcal{E}_{[M]}=\mathcal{E}_{[M^{\on{lc}}]}$.
\qed\enddemo

\begin{lemma}\label{nonquasiremarkslemma}
Let $(a_p)_{p\ge 1}$ be a decreasing sequence of positive real numbers with $\sum_{p\ge 1}a_p<+\infty$. Then $pa_p\rightarrow0$ as $p\rightarrow\infty$.
\end{lemma}


\subsection{The general case $\mathcal{E}_{[\mathcal{M}]}$}\label{subsection42}
Proposition \ref{nonquasiremarks} shows that $\mathcal{E}_{\{M\}}$ is non-quasianalytic if and only if $\mathcal{E}_{(M)}$ is. In the general case this is not true, e.g. let $\mathcal{M}=\{M^1,M^2\}$ such that $M^1\le M^2$, $\mathcal{E}_{[M^1]}$ is quasianalytic whereas $\mathcal{E}_{[M^2]}$ is not (take $M^1_p:=p!$ and $M^2_p:=p!^s$ for some $s>1$).

We prove now Theorem \ref{Matrix-non-quasi-analyticity}. The Roumieu part is obvious and the Beurling part will follow from the following Proposition \ref{Beurling-non-quasi-analyticity} which uses the idea of \cite[Lemma 5.1]{intersectionpaper}. We construct a non-quasianalytic sequence $N$ which is smaller than any sequence in the matrix $\mathcal{M}$. More precisely, we will show that $N\hyperlink{Mtriangle}{\vartriangleleft}\mathcal{M}$, while in \cite[Lemma 5.1]{intersectionpaper} only $N\hyperlink{Mbeurpreceq}{(\preceq)}\mathcal{M}$ was proved. Moreover the assumptions in \cite{intersectionpaper} where each $M^x$ is log-convex and $\mu^y_p\le\mu^x_p$ for all $p\in\NN$ and $y\le x$ will be not needed for our proof.

\begin{proposition}\label{Beurling-non-quasi-analyticity}
Let $\mathcal{M}:=\{M^x\in\RR_{>0}^{\NN}: x\in\Lambda=\RR_{>0}\}$ satisfy \hyperlink{Marb}{$(\mathcal{M})$} such that $\mathcal{E}_{[M^x]}$ is non-quasianalytic for each $x>0$. Then we get:
\begin{itemize}
\item[$(i)$] There exists $N$ with $N_0=1$ and $N^I=N$, $\mathcal{E}_{[N]}$ is non-quasianalytic and $N\hyperlink{Mtriangle}{\vartriangleleft}\mathcal{M}$, so $\mathcal{E}_{\{N\}}\subseteq\mathcal{E}_{(\mathcal{M})}$.

\item[$(ii)$] Let $U$ be a non-empty open subset of $\RR^r$. For every bounded subset $B$ in $\mathcal{E}_{(\mathcal{M})}(U)$ there exists a sequence $N$ as in $(i)$ such that $B$ is a bounded subset in $\mathcal{E}_{\{N\}}(U)$, too.

\item[$(iii)$] Let $\mathcal{N}:=\{N^x:x\in\RR_{>0}\}$ satisfy \hyperlink{Marb}{$(\mathcal{M})$} and \hyperlink{R-Comega}{$(\mathcal{M}_{\{\mathcal{C}^{\omega}\}})$} and such that $\mathcal{N}\hyperlink{Mtriangle}{\vartriangleleft}\mathcal{M}$. Then there exists a sequence $L$ which satisfies \hyperlink{lc}{$(\on{lc})$}, \hyperlink{mnq}{$(\on{nq})$} and finally $\mathcal{N}\hyperlink{Mtriangle}{\vartriangleleft}L\hyperlink{Mtriangle}{\vartriangleleft}\mathcal{M}$.
\end{itemize}
\end{proposition}

\demo{Proof}
$(i)$ Since $\Lambda=\RR_{>0}$ and $\mathcal{M}$ satisfies \hyperlink{Marb}{$(\mathcal{M})$} we can restrict to $\Lambda=\{\frac{1}{n}: n\in\NN_{>0}\}$, see \ref{weighmatrixdefinition}. By \cite[Theorem 1.3.8]{hoermander} and \cite[Theorem 4.2]{Komatsu73} we get
$$\forall\;x\in\NN_{>0}:\;\sum_{p\ge 1}\frac{1}{((M^{1/x})^I_p)^{1/p}}<+\infty.$$
Now, as in \cite[Lemma 5.1]{intersectionpaper}, we introduce sequences $(a_q)_{q\ge 0}$ and $(b_q)_{q\ge 0}$ recursively as follows. Put $a_0=b_0=0$, then let $a_q$ be the first integer such that
\begin{equation}\label{intersectionequ1}
a_q>b_{q-1},\hspace{30pt}\sum_{p=a_q+1}^{\infty}\frac{1}{((M^{1/(q+1)})^I_p)^{1/p}}\le\frac{2^{-q}}{q+1}.
\end{equation}
$b_q$ shall be the first integer such that $\frac{1}{q}((M^{1/q})^I_{a_q})^{1/a_q}<\frac{1}{q+1}((M^{1/(q+1)})^I_{b_q})^{1/b_q}$ holds. Since for each $q\in\NN_{>0}$ separately $p\mapsto((M^{1/q})^I_p)^{1/p}$ is increasing, tending to infinity and since $((M^{1/q})^I_p)^{1/p}\ge((M^{1/(q+1)})^I_p)^{1/p}$ for each $p,q\ge 1$ we have $a_q<b_q$ for each $q$.

Now introduce $N=(N_p)_p$ as follows. We put $N_0:=1$ and for $p\in\NN_{>0}$ we set
$$(N_p)^{1/p}=\frac{1}{q}((M^{1/q})^I_p)^{1/p}\;\;\text{for}\;\;b_{q-1}\le p\le a_q,\hspace{25pt}(N_p)^{1/p}=\frac{1}{q}((M^{1/q})^I_{a_q})^{1/a_q}\;\;\text{for}\;\;a_q+1\le p\le b_q-1.$$
{\itshape Claim.} The mapping $p\mapsto(N_p)^{1/p}$ is increasing, i.e. $N^I=N$.

If $b_{q-1}\le p<a_q$ and $a_q+1\le p<b_q-1$, then $(N_p)^{1/p}\le(N_{p+1})^{1/(p+1)}$ holds by definition. If $p=a_q$, then $(N_p)^{1/p}=\frac{1}{q}((M^{1/q})^I_p)^{1/p}\le\frac{1}{q}((M^{1/q})^I_{a_q})^{1/a_q}=(N_{p+1})^{1/(p+1)}$ and if $p=b_q-1$, then $(N_p)^{1/p}=\frac{1}{q}((M^{1/q})^I_{a_q})^{1/a_q}\le\frac{1}{q+1}((M^{1/(q+1)})^I_{b_q})^{1/b_q}=(N_{p+1})^{1/(p+1)}$ holds by the choice of $(b_q)_q$.

{\itshape Claim.} $\mathcal{E}_{[N]}$ is non-quasianalytic. First we have
\begin{align*}
&\sum_{p=a_1+1}^{\infty}\frac{1}{(N_p)^{1/p}}=\sum_{q=1}^{\infty}\left(\sum_{p=a_q+1}^{b_q-1}\frac{1}{(N_p)^{1/p}}+\sum_{p=b_q}^{a_{q+1}}\frac{1}{(N_p)^{1/p}}\right)
\\&
=\sum_{q=1}^{\infty}\left(\sum_{p=a_q+1}^{b_q-1}\frac{q}{((M^{1/q})^I_{a_q})^{1/a_q}}+\sum_{p=b_q}^{a_{q+1}}\frac{q+1}{((M^{1/(q+1)})^I_p)^{1/p}}\right)\underbrace{\le}_{(\star)}\sum_{q=1}^{\infty}\underbrace{\sum_{p=a_q+1}^{\infty}\frac{q+1}{((M^{1/(q+1)})^I_p)^{1/p}}}_{\le 2^{-q}}\le 1.
\end{align*}
$(\star)$ holds because by the choice of $(b_q)_q$ we have $\frac{1}{q+1}((M^{1/(q+1)})^I_p)^{1/p}\le\frac{1}{q}((M^{1/q})^I_{a_q})^{1/a_q}$ for $a_q+1\le p\le b_q-1$. Since $N=N^I$ and by \cite[Theorem 1.3.8]{hoermander} and \cite[Theorem 4.2]{Komatsu73} we are done.

{\itshape Claim.} $N\hyperlink{Mtriangle}{\vartriangleleft}\mathcal{M}$, i.e. $N\hyperlink{mtriangle}{\vartriangleleft}M^{1/x}$ for all $x\in\NN_{>0}$.

We have $(N_p)^{1/p}\le\frac{1}{q}((M^{1/q})^I_{p})^{1/p}\le\frac{1}{q}(M^{1/q}_p)^{1/p}$ whenever $p\ge b_{q-1}$, so $\mathcal{E}_{\{N\}}\subseteq\mathcal{E}_{(\mathcal{M})}$ follows.\vspace{6pt}

$(ii)$ Let $(K_j)_{j\in\NN_{>0}}$ be a fundamental system of compact subsets of $U$. For $j\in\NN_{>0}$ put
$$k_j:=\sup_{f\in B,i\in\NN, x\in K_j}\frac{2^{2ji} j^i\|f^{(i)}(x)\|_{L^i(\RR^r,\RR)}}{(M^{1/j})^I_{i}}.$$
Now introduce $(a_q)_q$ and $(b_q)_q$ as in $(i)$ but such that $a_q$ is the first integer satisfying \eqref{intersectionequ1} and additionally $k_{q+1} 2^{-a_q}\le 1$.

Let $\|\cdot\|_{N,K,h}$ be any fundamental continuous semi-norm in $\mathcal{E}_{(N)}$, then there exists $k\in\NN$ with $h^{-1}<2^k$ and $K\subseteq K_k$. For all $i\in\NN$ with $i>a_k$ there exists a unique $j>k$ with $a_{j-1}<i\le a_j$. By definition this implies $\frac{1}{j}((M^{1/j})^I_p)^{1/p}\le(N_p)^{1/p}$ for all $p\in\NN$ with $p\le i$ and so $(M^{1/j})^I_p\le j^p N_p$ for such $p$. Thus we get for all $i$ sufficiently large
\begin{align*}
\sup_{x\in K}\frac{\|f^{(i)}(x)\|_{L^i(\RR^r,\RR)}}{h^i N_i}\le\sup_{x\in K_j}\frac{2^{k i} j^i\|f^{(i)}(x)\|_{L^i(\RR^r,\RR)}}{(M^{1/j})^I_i}\le 2^{-ji}k_j\le 2^{-a_{j-1}} k_j\le 1,
\end{align*}
for all $f\in B$.

We are done since by Proposition \ref{nonquasiremarks} the matrix $\mathcal{M}$ has \hyperlink{B-Comega}{$(\mathcal{M}_{(\mathcal{C}^{\omega})})$} and so for each $M^{1/j}$ separately we get $\mathcal{E}_{[(M^{1/j})^{\on{lc}}]}=\mathcal{E}_{[(M^{1/j})^I]}=\mathcal{E}_{[M^{1/j}]}$.\vspace{6pt}

$(iii)$ By \hyperlink{R-Comega}{$(\mathcal{M}_{\{\mathcal{C}^{\omega}\}})$} and \cite[Theorem 2.15 (1)]{compositionpaper} we can assume that each $N^x\in\mathcal{N}$ is log-convex since $\mathcal{E}_{\{N^x\}}=\mathcal{E}_{\{(N^x)^{\on{lc}}\}}$ for all $x$ (we can drop all small indices for which possibly $\liminf_{k\rightarrow\infty}(n^x_k)^{1/k}=0$ without changing the space $\mathcal{E}_{\{\mathcal{N}\}}$).

By Proposition \ref{nonquasiremarks} and $(i)$ there exists $P$ with $P\hyperlink{Mtriangle}{\vartriangleleft}\mathcal{M}$, $\mathcal{E}_{[P]}=\mathcal{E}_{[P^{\on{lc}}]}$ and $P^{\on{lc}}$ has \hyperlink{mnq}{$(\on{nq})$}. Consequently $\mathcal{E}_{[P^{\on{lc}}]}$ is non-quasianalytic and $P^{\on{lc}}\hyperlink{Mtriangle}{\vartriangleleft}\mathcal{M}$ holds, too.

On the other hand by \cite[Lemma 3.5.7]{dissertationesser} there exists $Q$ with $\mathcal{N}\hyperlink{Mtriangle}{\vartriangleleft}Q\hyperlink{Mtriangle}{\vartriangleleft}\mathcal{M}$.

Now put $Q'_k:=\max\{P^{\on{lc}}_k,Q_k\}$. Since $Q'\ge P^{\on{lc}}$ we have that $\mathcal{E}_{[Q']}$ is non-quasianalytic, $\mathcal{E}_{[Q']}=\mathcal{E}_{[Q'^{\on{lc}}]}$ and $Q'^{\on{lc}}$ satisfies \hyperlink{mnq}{$(\on{nq})$}.

On the other hand $Q'\ge Q$ implies $\mathcal{N}\hyperlink{Mtriangle}{\vartriangleleft}Q'$. Since $\mathcal{E}_{[Q']}=\mathcal{E}_{[Q'^{\on{lc}}]}$ and each $N^x\in\mathcal{N}$ has \hyperlink{lc}{$(\on{lc})$}, also $\mathcal{N}\hyperlink{Mtriangle}{\vartriangleleft}Q'^{\on{lc}}$ follows.

Finally $Q'^{\on{lc}}\hyperlink{Mtriangle}{\vartriangleleft}\mathcal{M}$ holds because $Q'^{\on{lc}}\le Q'$ and $P^{\on{lc}},Q\hyperlink{Mtriangle}{\vartriangleleft}\mathcal{M}$.

The conclusion follows now by defining $L:=Q'^{\on{lc}}$.
\qed\enddemo

If $\mathcal{M}=\Omega$ is coming from $\omega\in\hyperlink{omset1}{\mathcal{W}}$, then we obtain the following consequence:

\begin{corollary}\label{nonquasiweight}
Let $\omega\in\hyperlink{omset1}{\mathcal{W}}$ be given, TFAE:

\begin{itemize}
\item[$(i)$] $\omega$ has \hyperlink{omnq}{$(\omega_{\on{nq}})$},

\item[$(ii)$] $\mathcal{E}_{\{\omega\}}$ contains functions with compact support,

\item[$(iii)$] $\mathcal{E}_{(\omega)}$ contains functions with compact support,

\item[$(iv)$] some $\Omega^l$ has \hyperlink{mnq}{$(\on{nq})$},

\item[$(v)$] each $\Omega^l$ has \hyperlink{mnq}{$(\on{nq})$}.
\end{itemize}
\end{corollary}

\demo{Proof}
By \cite[5.5]{compositionpaper} the matrix $\Omega$ is \hyperlink{Msc}{$(\mathcal{M}_{\on{sc}})$}. By \cite[5.5, Corollary 5.8 $(1)$]{compositionpaper} we have $(i)\Leftrightarrow(iv)\Leftrightarrow(v)$. The rest follows from Theorem \ref{Matrix-non-quasi-analyticity}.
\qed\enddemo

\section{Characterization of $\mathcal{E}_{[\mathcal{M}]}$ using the Fourier transform}\label{section5}
Using the central results from Sections \ref{section3} and \ref{section4} we are now able to characterize functions in $\mathcal{E}_{[\mathcal{M}]}$ in terms of the decay of its Fourier transform. First put
$$\mathcal{D}(\RR^r):=\{f\in\mathcal{E}(\RR^r):\;\exists\;K\subset\subset\RR^r,\;\supp(f)\subseteq K\}.$$
Let $\mathcal{M}=\{M^x: x\in\Lambda\}$ satisfy \hyperlink{Marb}{$(\mathcal{M})$}. If $\mathcal{E}_{(\mathcal{M})}$ respectively $\mathcal{E}_{\{\mathcal{M}\}}$ is non-quasianalytic, then
$$\mathcal{D}_{(\mathcal{M})}(U):=\{f\in\mathcal{E}(\RR^r):\;\exists\;K\subset\subset U\;\supp(f)\subseteq K,\;\forall\;x\in\Lambda\;\forall\;h>0\;:\|f\|_{M^x,\RR,h}<+\infty\}$$
respectively
$$\mathcal{D}_{\{\mathcal{M}\}}(K):=\{f\in\mathcal{E}(\RR^r):\;\exists\;K\subset\subset U\;\supp(f)\subseteq K,\exists\;x\in\Lambda\;\exists\;h>0\;:\|f\|_{M^x,\RR,h}<+\infty\}$$
is non-trivial.

On the other hand let $\mathcal{M}=\{M^x: x\in\Lambda\}$ be \hyperlink{Msc}{$(\mathcal{M}_{\on{sc}})$} and let $K\subset\subset\RR^r$ be compact. Then for $x\in\Lambda$ and $h>0$ introduce the Banach space
$$\hat{\mathcal{D}}_{x,h}(K):=\{f\in\mathcal{E}(\RR^r):\;\supp(f)\subseteq K, \|f\|^{\hat{}}_{x,h}<+\infty\},$$
where $\|f\|^{\hat{}}_{x,h}:=\int_{\RR^r}|\hat{f}(t)|\exp(h\omega_{M^x}(t))dt$. So one can define
$$\hat{\mathcal{D}}_{(\omega_{\mathcal{M}})}(K):=\underset{x\in\Lambda, h>0}{\varprojlim}\hat{\mathcal{D}}_{x,h}(K)\hspace{30pt}\hat{\mathcal{D}}_{\{\omega_{\mathcal{M}}\}}(K):=\underset{x\in\Lambda, h>0}{\varinjlim}\hat{\mathcal{D}}_{x,h}(K),$$
and for non-empty open $U\subseteq\RR^r$
$$\hat{\mathcal{D}}_{(\omega_{\mathcal{M}})}(U):=\underset{K\subset\subset U}{\varinjlim}\hat{\mathcal{D}}_{(\omega_{\mathcal{M}})}(K)\hspace{30pt}\hat{\mathcal{D}}_{\{\omega_{\mathcal{M}}\}}(U):=\underset{K\subset\subset U}{\varinjlim}\hat{\mathcal{D}}_{\{\omega_{\mathcal{M}}\}}(K).$$
Now we formulate our main theorem:
\begin{theorem}\label{centralfouriertheorem}
Let $\mathcal{M}:=\{M^x: x\in\Lambda\}$ be \hyperlink{Msc}{$(\mathcal{M}_{\on{sc}})$}. Moreover assume that
\begin{itemize}
\item[$(i)$] $\mathcal{M}$ has $(\mathcal{M}_{[\on{L}]})$,
\item[$(ii)$] $\mathcal{M}$ has $(\mathcal{M}_{[\on{mg}]})$,
\item[$(iii)$] $\mathcal{E}_{[\mathcal{M}]}$ is non-quasianalytic.
\end{itemize}
Then we obtain the equalities
$$\mathcal{D}_{[\mathcal{M}]}=\mathcal{D}_{[\omega_{\mathcal{M}}]}=\hat{\mathcal{D}}_{[\omega_{\mathcal{M}}]}.$$
\end{theorem}

{\itshape Examples.} The previous theorem is valid if $\mathcal{M}=\Omega$ for some $\omega\in\hyperlink{omset1}{\mathcal{W}}$ with \hyperlink{omnq}{$(\omega_{\on{nq}})$} or also for the Gevrey-matrix $\mathcal{G}$.

For the proof we have to generalize \cite[Lemma 3.3]{BraunMeiseTaylor90}. Let $K\subset\subset\RR^r$ and let $H_K(t):=\sup_{s\in K}\langle t,s\rangle$ be the support function. $\lambda_r(K)$ shall denote the Lebesgue measure of $K$.

\begin{lemma}\label{centralfourierlemma}
Let $\mathcal{M}=\{M^x: x\in\Lambda\}$ be \hyperlink{Msc}{$(\mathcal{M}_{\on{sc}})$} and $f\in\mathcal{D}(\RR^r)$.

\begin{itemize}
\item[$(i)$] Let $x\in\Lambda$ and $h>0$ be arbitrary and assume that $\|f\|^{\hat{}}_{x,h}=:C<+\infty$. Then
    \begin{equation}\label{bmtequ1}
    \sup_{\alpha\in\NN^r, t\in\RR^r}|f^{(\alpha)}(t)|\exp\left(-h\varphi^{*}_{\omega_{M^x}}\left(\frac{|\alpha|}{h}\right)\right)\le\frac{C}{(2\pi)^r}.
    \end{equation}
    holds.
\item[$(ii)$] Let $\mathcal{M}$ satisfy additionally $(\mathcal{M}_{[\on{L}]})$.

In the Roumieu case assume that there exist some $x\in\Lambda$ and $C,h>0$ such that \eqref{bmtequ1} is valid. Then there exists $D\ge 1$ depending on $x,h$ and the dimension $r$ and there exist $y\in\Lambda$ and $L\ge 1$ depending only on $x$ and $r$ such that with $K:=\supp(f)$ we have for all $z\in\CC^r$
\begin{equation}\label{bmtequ2}
|\hat{f}(z)|\le\lambda_r(K)\frac{CD}{(2\pi)^r}\exp\left(H_K(\Im(z))-\frac{h}{L}\omega_{M^y}(z)\right).
\end{equation}
\end{itemize}
In the Beurling case for arbitrary $y\in\Lambda$ and $h>0$ there exists $D\ge 1$ depending on $x,h$ and the dimension $r$ and there exist $x\in\Lambda$ and $L\ge 1$ depending only on $y$ and $r$ such that \eqref{bmtequ2} holds (with $y,D,L$) provided \eqref{bmtequ1} is valid (with $x,h,C$).
\end{lemma}
For $(ii)$ it is sufficient to assume \eqref{R-L-consequ} in the Roumieu and \eqref{B-L-consequ} in the Beurling case, see Proposition \ref{2}.

\demo{Proof}
$(i)$ Since each $\omega_{M^x}\in\hyperlink{omset0}{\mathcal{W}_0}$ we can replace in the proof of \cite[Lemma 3.3 $(1)$]{BraunMeiseTaylor90} the weight $\omega$ by $\omega_{M^x}$.

$(ii)$ We consider the Roumieu case. Iterating \eqref{R-L-consequ} yields $\omega_{M^y}(rt)\le\frac{L}{2}\omega_{M^x}(t)+\frac{L}{2}$ for all $t\ge 0$ and for some $y\in\Lambda$ and $L\ge 1$ both depending only on $x$ and $r$. By \hyperlink{om3}{$(\omega_3)$} for $\omega_{M^y}$ there exists some $B\ge 1$ such that $(2h/L)\omega_{M^y}(t)-\log(t)\ge(h/L)\omega_{M^y}(t)-B$ for all $t\ge 1$.

Then follow \cite[Lemma 3.3 $(2)$]{BraunMeiseTaylor90}.
\qed\enddemo

Lemma \ref{centralfourierlemma} and the Paley-Wiener theorem for $\mathcal{D}(K)$ (see \cite[7.3.1]{hoermander}) imply

\begin{proposition}\label{centralfourierproposition}
Let $\mathcal{M}=\{M^x: x\in\Lambda\}$ be \hyperlink{Msc}{$(\mathcal{M}_{\on{sc}})$} with $(\mathcal{M}_{[\on{L}]})$, let $K\subset\subset\RR^r$ be a compact convex set and $f\in L^1(\RR^r)$.
\begin{itemize}
\item[$(i)$] The Roumieu case. The following are equivalent:
\begin{itemize}
\item[$(a)$] $f\in\hat{\mathcal{D}}_{\{\omega_{\mathcal{M}}\}}(K)$,

\item[$(b)$] $f\in\mathcal{D}(K)$ and there exists $x\in\Lambda$ and $l>0$ such that $\|f\|_{\omega_{M^x},K,l}<+\infty$,

\item[$(c)$] there exist $x\in\Lambda$ and $C,l>0$ such that for all $z\in\CC^r$ we have
$$|\hat{f}(z)|\le C\exp(H_K(\Im(z))-l\omega_{M^x}(z)).$$
\end{itemize}

\item[$(ii)$] The Beurling case. The following are equivalent:
\begin{itemize}
\item[$(a)$] $f\in\hat{\mathcal{D}}_{(\omega_{\mathcal{M}})}(K)$,

\item[$(b)$] $f\in\mathcal{D}(K)$ and for all $x\in\Lambda$ and $l>0$ we have $\|f\|_{\omega_{M^x},K,l}<+\infty$,

\item[$(c)$] for all $x\in\Lambda$ and $l>0$ there exists $C\ge 1$ such that for all $z\in\CC^r$ we have
$$|\hat{f}(z)|\le C\exp(H_K(\Im(z))-l\omega_{M^x}(z)).$$
\end{itemize}
\end{itemize}
\end{proposition}

Theorem \ref{centralfouriertheorem} follows now by applying Theorem \ref{4} and Proposition \ref{centralfourierproposition}.

\section{Comparison of the classes $\mathcal{E}_{[M]}$ and $\mathcal{E}_{[\omega]}$}\label{section6}
In \cite{BonetMeiseMelikhov07} the authors compared the classical methods which are used to introduce classes of ultradifferentiable functions, either by a weight sequence $M$ or a weight function $\omega$. In \cite{compositionpaper} we have introduced the technique of associating a weight matrix $\Omega$ to a given function $\omega$. The aim of this section is to reformulate the comparison theorems in view of this new method.

\begin{theorem}\label{finalcomparisontheorem}
Let $\omega\in\hyperlink{omset1}{\mathcal{W}}$, TFAE:

\begin{itemize}
\item[$(i)$] There exists $N\in\hyperlink{LCset}{\mathcal{LC}}$ with $\mathcal{E}_{[N]}=\mathcal{E}_{[\omega]}=\mathcal{E}_{[\Omega]}$,

\item[$(ii)$] $\omega$ has \hyperlink{om6}{$(\omega_6)$},

\item[$(iii)$] there exists $N\in\hyperlink{LCset}{\mathcal{LC}}$ such that for each $l>0$ we have $\mathcal{E}_{[\Omega^l]}=\mathcal{E}_{[N]}$ or equivalently $N\hyperlink{approx}{\approx}\Omega^l$.
\end{itemize}

Additionally we have:

\begin{itemize}
\item[$(a)$] $N$ and each $\Omega^l$ satisfy \hyperlink{mg}{$(\on{mg})$}.

\item[$(b)$] If $\omega$ has \hyperlink{om2}{$(\omega_2)$}, then $\liminf_{p\rightarrow\infty}(n_p)^{1/p}>0$, \hyperlink{holom}{$(\mathcal{M}_{\mathcal{H}})$} for $\Omega$ and
\begin{itemize}
\item[$*$] $\omega_N$ and each $\omega_{\Omega^l}$ satisfy \hyperlink{om2}{$(\omega_2)$},
\item[$*$] $N$ and each $\Omega^l$ have \hyperlink{beta3}{$(\beta_3)$},
\item[$*$] $\mathcal{E}_{[\omega_N]}=\mathcal{E}_{[N]}=\mathcal{E}_{[\omega]}=\mathcal{E}_{[\Omega^l]}$ for each $l>0$.
\end{itemize}
If $\omega$ has \hyperlink{om5}{$(\omega_5)$}, then $\lim_{p\rightarrow\infty}(n_p)^{1/p}=\infty$, \hyperlink{B-Comega}{$(\mathcal{M}_{(\mathcal{C}^{\omega})})$} for $\Omega$ and $\omega_N$ and each $\omega_{\Omega^l}$ satisfy \hyperlink{om5}{$(\omega_5)$}.
\end{itemize}
\end{theorem}

In the next theorem we start with a weight sequence $N$ and not with a weight function $\omega$ as before.

\begin{theorem}\label{finalcomparisontheorem1}
Let $N\in\hyperlink{LCset}{\mathcal{LC}}$ with \hyperlink{beta3}{$(\beta_3)$}, TFAE:

\begin{itemize}
\item[$(i)$] There exists $\omega\in\hyperlink{omset1}{\mathcal{W}}$ such that $\mathcal{E}_{[\omega]}=\mathcal{E}_{[N]}$,

\item[$(ii)$] $N$ satisfies \hyperlink{mg}{$(\on{mg})$},

\item[$(iii)$] $\mathcal{E}_{[\omega_N]}=\mathcal{E}_{[N]}$ holds.
\end{itemize}

Let $\Omega:=\{\Omega^l: l>0\}$ be the matrix associated to $\omega$ arising in $(i)$. We get for each $l>0$:

\begin{itemize}
\item[$(a)$] $\omega,\omega_{\Omega^l},\omega_N\in\hyperlink{omset1}{\mathcal{W}}$ satisfy \hyperlink{om6}{$(\omega_6)$},

\item[$(b)$] $\omega_{\Omega^l}\hyperlink{sim}{\sim}\omega\hyperlink{sim}{\sim}\omega_N$,

\item[$(c)$] $\mathcal{E}_{[\omega_N]}=\mathcal{E}_{[N]}=\mathcal{E}_{[\omega]}=\mathcal{E}_{[\Omega^l]}$,

\item[$(d)$] $N\hyperlink{approx}{\approx}\Omega^l$,

\item[$(e)$] $\Omega^l$ has \hyperlink{mg}{$(\on{mg})$}.

\item[$(f)$] If $N$ satisfies $\liminf_{p\rightarrow\infty}(n_p)^{1/p}>0$, then
\begin{itemize}
\item[$(*)$] $\omega,\omega_N$ and each $\omega_{\Omega^l}$ have \hyperlink{om2}{$(\omega_2)$},
\item[$(*)$] each $\Omega^l$ has \hyperlink{beta3}{$(\beta_3)$},
\item[$(*)$] \hyperlink{holom}{$(\mathcal{M}_{\mathcal{H}})$} for $\Omega$.
\end{itemize}
If $N$ satisfies $\lim_{p\rightarrow\infty}(n_p)^{1/p}=+\infty$, then
\begin{itemize}
\item[$(*)$] $\omega,\omega_N$ and each $\omega_{\Omega^l}$ have \hyperlink{om5}{$(\omega_5)$},
\item[$(*)$] \hyperlink{B-Comega}{$(\mathcal{M}_{(\mathcal{C}^{\omega})})$} for $\Omega$.
\end{itemize}
\end{itemize}
\end{theorem}
Theorem \ref{finalcomparisontheorem} and Theorem \ref{finalcomparisontheorem1} follow by the results below, \cite[Section 5]{compositionpaper} and \cite{BonetMeiseMelikhov07}, see also \cite[6.1-6.4]{dissertation}.

\begin{theorem}\label{comparisontheorem}
Let $\omega\in\hyperlink{omset1}{\mathcal{W}}$, $U\subseteq\RR^r$ non-empty open. Then we get:
\begin{itemize}
\item[$(1)$] $\omega$ has \hyperlink{om6}{$(\omega_6)$} if and only if $\mathcal{E}_{[\Omega^l]}(U)=\mathcal{E}_{[\omega]}(U)$ holds for each $l>0$. Moreover for each $l>0$
\begin{itemize}
\item[$(a)$] $\omega\hyperlink{sim}{\sim}\omega_{\Omega^l}$,
\item[$(b)$] $\Omega^l\in\hyperlink{LCset}{\mathcal{LC}}$,
\item[$(c)$] $\omega_{\Omega^l}\in\hyperlink{omset1}{\mathcal{W}}$ with \hyperlink{om6}{$(\omega_6)$},
\item[$(d)$] $\Omega^x\hyperlink{approx}{\approx}\Omega^y$ holds for all $x,y>0$,
\item[$(e)$] $\Omega^l$ satisfies \hyperlink{mg}{$(\on{mg})$}.
\end{itemize}
\item[$(2)$] Let $\omega$ be as in $(1)$ with \hyperlink{om2}{$(\omega_2)$}, then
\begin{itemize}
\item[$(a)$] $\Omega$ has \hyperlink{holom}{$(\mathcal{M}_{\mathcal{H}})$},
\item[$(b)$] each $\Omega^l$ satisfies \hyperlink{beta3}{$(\beta_3)$},
\item[$(c)$] each $\omega_{\Omega^l}$ has \hyperlink{om2}{$(\omega_2)$}.
\end{itemize}
If $\omega$ is as in $(1)$ with \hyperlink{om5}{$(\omega_5)$}, then
\begin{itemize}
\item[$(d)$] $\Omega$ has \hyperlink{B-Comega}{$(\mathcal{M}_{(\mathcal{C}^{\omega})})$},
\item[$(e)$] each $\omega_{\Omega^l}$ has \hyperlink{om5}{$(\omega_5)$}.
\end{itemize}
\end{itemize}
\end{theorem}

\demo{Proof}
$(1)$ This was already shown in \cite[Section 5]{compositionpaper}.

$(2)$ To prove \hyperlink{beta3}{$(\beta_3)$} for each $\Omega^l$ we proceed similarly as in \cite[Lemma 12  $(1)\Rightarrow(2)$]{BonetMeiseMelikhov07} (each sequence $\Omega^l$ satisfies the required assumptions).

If $\omega$ has \hyperlink{om2}{$(\omega_2)$} or \hyperlink{om5}{$(\omega_5)$}, then by \cite[Lemma 5.7]{compositionpaper} each $\omega_{\Omega^l}$ too and by \cite[Proposition 4.6 $(1)$, Corollary 5.15]{compositionpaper} we get \hyperlink{holom}{$(\mathcal{M}_{\mathcal{H}})$} or \hyperlink{B-Comega}{$(\mathcal{M}_{(\mathcal{C}^{\omega})})$} for $\Omega$.
\qed\enddemo

In the next result we start with a weight sequence $M$ and not with $\omega$.

\begin{theorem}\label{conversecomparison}
Let $M\in\hyperlink{LCset}{\mathcal{LC}}$ with \hyperlink{beta3}{$(\beta_3)$} and \hyperlink{mg}{$(\on{mg})$}. Let $r\in\NN_{>0}$ and $U\subseteq\RR^r$ be non-empty open. Then
\begin{itemize}
\item[$(1)$] $\omega_{M}\in\hyperlink{omset1}{\mathcal{W}}$ has \hyperlink{om6}{$(\omega_6)$}.

\item[$(2)$] $\mathcal{E}_{[\omega_{M}]}(U)=\mathcal{E}_{[N^l]}(U)=\mathcal{E}_{[M]}(U)$ for each $l>0$, where $N^l_p:=\exp(\frac{1}{l}\varphi^{*}_{\omega_M}(lp))$. Moreover $N^1=M$ and for each $l>0$
\begin{itemize}
\item[$(a)$] $N^l\in\hyperlink{LCset}{\mathcal{LC}}$ and has \hyperlink{mg}{$(\on{mg})$},
\item[$(b)$] $\omega_{N^l}\hyperlink{sim}{\sim}\omega_M$, $\omega_{N^l}\in\hyperlink{omset1}{\mathcal{W}}$ with \hyperlink{om6}{$(\omega_6)$},
\item[$(c)$] $M\hyperlink{approx}{\approx}N^l$.
\end{itemize}
\item[$(3)$] If $M$ satisfies $\liminf_{p\rightarrow\infty}(m_p)^{1/p}>0$, then
\begin{itemize}
\item[$(a)$] \hyperlink{om2}{$(\omega_2)$} for $\omega_M$ and each $\omega_{N^l}$,
\item[$(b)$] each $N^l$ has \hyperlink{beta3}{$(\beta_3)$} and $\liminf_{p\rightarrow\infty}(n^l_p)^{1/p}>0$.
\end{itemize}
If $M$ satisfies $\lim_{p\rightarrow\infty}(m_p)^{1/p}=\infty$, then
\begin{itemize}
\item[$(c)$] \hyperlink{om5}{$(\omega_5)$} for $\omega_M$ and each $\omega_{N^l}$,
\item[$(d)$] each $N^l$ has $\lim_{p\rightarrow\infty}(n^l_p)^{1/p}=+\infty$.
\end{itemize}
\end{itemize}
\end{theorem}

\demo{Proof}
$(1)$ By \ref{assofuncproper} we get $\omega_{M}\in\hyperlink{omset0}{\mathcal{W}_0}$, by \cite[Lemma 12  $(2)\Rightarrow(4)$]{BonetMeiseMelikhov07} we get \hyperlink{om1}{$(\omega_1)$} and by \cite[Proposition 3.6]{Komatsu73} we get \hyperlink{om6}{$(\omega_6)$} for $\omega_M$.\vspace{6pt}

$(2)$ In Theorem \ref{comparisontheorem} consider $\omega=\omega_M$ and then $\mathcal{E}_{[\omega_{M}]}(U)=\mathcal{E}_{[N^l]}(U)$ for each $l>0$. By \cite[5.5]{compositionpaper} we have $N^l\in\hyperlink{LCset}{\mathcal{LC}}$ and so
\begin{align*}
M_p&=\sup_{t\ge0}\frac{t^p}{\exp(\omega_M(t))}=\exp\left(\sup_{t\ge0}\left(p\log(t)-\omega_M(t)\right)\right)=\exp\left(\varphi^{*}_{\omega_M}(p)\right)=:N^1_p,
\end{align*}
for all $p\in\NN$. Thus $\mathcal{E}_{[M]}=\mathcal{E}_{[N^1]}=\mathcal{E}_{[\omega_M]}=\mathcal{E}_{[N^l]}$ which implies $M\hyperlink{approx}{\approx}N^l$ and \hyperlink{mg}{$(\text{mg})$} follows for each $N^l$.

By \ref{assofuncproper} we have $\omega_{N^l}\in\hyperlink{omset0}{\mathcal{W}_0}$, hence \cite[Lemma 5.7]{compositionpaper} applied to $\omega_M$ implies $\omega_{N^l}\hyperlink{sim}{\sim}\omega_M$ for each $l>0$ and so \hyperlink{om1}{$(\omega_1)$} and \hyperlink{om6}{$(\omega_6)$} for each $\omega_{N^l}$ follow.\vspace{6pt}

$(3)$ By \ref{assofuncproper} the assumption $\liminf_{p\rightarrow\infty}(m_p)^{1/p}>0$ implies \hyperlink{om2}{$(\omega_2)$} for $\omega_M$ and each $\omega_{N^l}$. Again by \cite[Proposition 4.6 $(1)$, Corollary 5.15]{compositionpaper} we get $\liminf_{p\rightarrow\infty}(n^l_p)^{1/p}>0$ for each $l>0$ and similarly for $\lim_{p\rightarrow\infty}(m_p)^{1/p}=+\infty$.

To show \hyperlink{beta3}{$(\beta_3)$} for each $N^l$ we follow again \cite[Lemma 12 $(1)\Rightarrow(2)$]{BonetMeiseMelikhov07}.
\qed\enddemo

\appendix
\section{Nuclearity of the connecting mappings for $\mathcal{E}_{[\mathcal{M}]}$}
First we recall \cite[Lemma 2.3]{Komatsu73}:

\begin{lemma}\label{nuclearlemma}
The identity mapping
$$\mathcal{C}^{r+1}(K,\RR)\longrightarrow\mathcal{C}(K,\RR)$$
is nuclear for each compact set $K\subset\subset\RR^r$ with smooth boundary.
\end{lemma}

Let $\mathcal{M}:=\{M^x: x\in\Lambda\}$ be \hyperlink{Marb}{$(\mathcal{M})$}. For $x\le y$, $h\le k$ and a compact set $K\subset\subset\RR^r$ with smooth boundary consider the inclusion
\begin{equation}\label{nuclearinclusion}
\mathcal{E}_{M^x,h}(K,\RR)\longrightarrow\mathcal{E}_{M^y,k}(K,\RR),
\end{equation}
and we are going to prove the matrix generalization of \cite[Proposition 2.4]{Komatsu73}:

\begin{proposition}\label{Nuclearproposition}
Let $\mathcal{M}$ satisfy \hyperlink{Marb}{$(\mathcal{M})$}.
\begin{itemize}
\item[$(a)$] If \hyperlink{R-dc}{$(\mathcal{M}_{\{\on{dc}\}})$}, then $\forall\;x\in\Lambda\;\forall\;h>0\;\exists\;y\in\Lambda\;\exists\;k>0:$ \eqref{nuclearinclusion} is nuclear.

\item[$(b)$] If \hyperlink{B-dc}{$(\mathcal{M}_{(\on{dc})})$}, then $\forall\;y\in\Lambda\;\forall\;k>0\;\exists\;x\in\Lambda\;\exists\;h>0:$ \eqref{nuclearinclusion} is nuclear.
\end{itemize}
\end{proposition}

\demo{Proof}
As already pointed out in \cite[Proposition 2.4]{Komatsu73}, since each inclusion mapping is a product of two inclusion mappings of the same type, it is enough to show quasi-nuclearity, see \cite[Theorem 3.3.2]{pietsch}. For convenience put $X:=\mathcal{E}_{M^x,h}(K,\RR)$ and $Y:=\mathcal{E}_{M^y,k}(K,\RR)$. So we have to show that there exists $(u_j)_j$, $u_j\in X'$, such that $\sum_{j=1}^{\infty}\|u_j\|_{X'}<+\infty$ and
$$\|f\|_Y\le\sum_{j=1}^{\infty}|\langle f,u_j\rangle_X|\hspace{20pt}\forall\;f\in X.$$
Now we point out that
\begin{equation}\label{nuclear1}
\|f\|_Y:=\sup_{\alpha\in\NN^r,x\in K}\frac{\left|f^{(\alpha)}(x)\right|}{k^{|\alpha|} M^y_{|\alpha|}}=\sup_{\alpha\in\NN^r}\frac{\|f^{(\alpha)}\|_{\mathcal{C}(K,\RR)}}{k^{|\alpha|} M^y_{|\alpha|}}\le\sum_{\alpha\in\NN^r}\frac{\|f^{(\alpha)}\|_{\mathcal{C}(K,\RR)}}{k^{|\alpha|} M^y_{|\alpha|}}.
\end{equation}
By Lemma \ref{nuclearlemma} there exists $(v_j)_j$, $v_j\in(\mathcal{C}^{r+1}(K,\RR))'$ such that
\begin{equation}\label{nuclear2}
\sum_{j=1}^{\infty}\|v_j\|_{(\mathcal{C}^{r+1}(K,\RR))'}<+\infty,\hspace{30pt}\|f^{(\alpha)}\|_{\mathcal{C}(K,\RR)}\le\sum_{j=1}^{\infty}\left|\langle f^{(\alpha)},v_j\rangle_{\mathcal{C}^{r+1}(K,\RR)}\right|.
\end{equation}
Now let $u_{\alpha,j}$ be the linear functional on $X$ defined by
\begin{equation}\label{nuclear3}
\langle f,u_{\alpha,j}\rangle:=\frac{\langle f^{(\alpha)},v_j\rangle_{\mathcal{C}^{r+1}(K,\RR)}}{k^{|\alpha|} M^y_{|\alpha|}}.
\end{equation}
By \eqref{nuclear1} and \eqref{nuclear2} we get:
$$\|f\|_Y\le\sum_{\alpha\in\NN^r,j\in\NN}\left|\langle f,u_{\alpha,j}\rangle\right|.$$
Moreover, by \eqref{nuclear3} we have
\begin{align*}
\left|\langle f,u_{\alpha,j}\rangle\right|&=\frac{|\langle f^{(\alpha)},v_j\rangle_{\mathcal{C}^{r+1}(K,\RR)}|}{k^{|\alpha|} M^y_{|\alpha|}}\le\frac{\|f^{(\alpha)}\|_{\mathcal{C}^{r+1}(K,\RR)}\|v_j\|_{\mathcal{C}^{r+1}(K,\RR)}}{k^{|\alpha|} M^y_{|\alpha|}}
\\&
\le\sup_{0\le|q|\le r+1}\frac{\|f^{(\alpha+q)}\|_{\mathcal{C}(K,\RR)}\|v_j\|_{\mathcal{C}^{r+1}(K,\RR)}}{k^{|\alpha|} M^y_{|\alpha|}}
\\&
\le\sup_{0\le|q|\le r+1}\frac{\|f\|_{X}h^{|\alpha+q|}M^x_{|\alpha+q|}\|v_j\|_{\mathcal{C}^{r+1}(K,\RR)}}{k^{|\alpha|} M^y_{|\alpha|}}
\\&
\le\frac{h^{|\alpha|}}{k^{|\alpha|}}\sup_{0\le|q|\le r+1}\frac{h^{|\alpha+q|} M^x_{|\alpha+q|}}{h^{|\alpha|} M^y_{|\alpha|}}\|f\|_X\|v_j\|_{\mathcal{C}^{r+1}(K,\RR)}.
\end{align*}
$(a)$ {\itshape Roumieu case.} By \hyperlink{R-dc}{$(\mathcal{M}_{\{\on{dc}\}})$} for given $x\in\Lambda$ we can find $x_1\in\Lambda$ and $H\ge 1$ such that $M^x_{|\alpha+q|}=M^x_{|\alpha|+|q|}\le H^{|\alpha|} M^{x_1}_{|\alpha|}$ for all $\alpha\in\NN^r$ and $q\in\NN^r$ with $0\le|q|\le r+1$. $M^y\ge M^{x_1}$ holds for $y\ge x_1$ and so
$$\sup_{0\le|q|\le r+1}\frac{h^{|\alpha+q|} M^x_{|\alpha+q|}}{h^{|\alpha|} M^y_{|\alpha|}}\le A H^{|\alpha|}(1+h^{r+1})$$
for some constant $A>0$. Hence if we choose $k$ such that $k>Hh\Leftrightarrow\frac{Hh}{k}<1$, then by \eqref{nuclear2} we get
$$\sum_{\alpha\in\NN^r, j\in\NN}\|u_{\alpha,j}\|_{X'}\le A\sum_{\alpha\in\NN^r, j\in\NN}\left(\frac{Hh}{k}\right)^{|\alpha|}\|v_j\|_{(\mathcal{C}^{r+1}(K,\RR))')}(1+h^{r+1})<+\infty.$$

$(b)$ {\itshape Beurling case.} By \hyperlink{B-dc}{$(\mathcal{M}_{(\on{dc})})$} for given $y\in\Lambda$ we can find $y_1\in\Lambda$ and $H\ge 1$ such that $M^{y_1}_{|\alpha+q|}\le H^{|\alpha|}M^{y}_{|\alpha|}$ for all $\alpha\in\NN^r$ and $q\in\RR^r$ with $0\le|q|\le r+1$.

So for given $y\in\Lambda$ and $k>0$ (both small) we can take $x\le y_1$, $h<\frac{k}{H}$ and estimate as for the Roumieu case.
\qed\enddemo


\begin{thebibliography}{99}

\bibitem{BonetMeiseMelikhov07} \textsc{J. Bonet, R. Meise, S. N. Melikhov}: \textit{A comparison of two different ways to define classes of ultradifferentiable
	functions}, Bull. Belg. Math. Soc. Simon Stevin, 2007, \textbf{14}, 424--444.

\bibitem{BraunMeiseTaylor90} \textsc{R. W. Braun, R. Meise, B. A. Taylor}: \textit{Ultradifferentiable functions and {F}ourier analysis}, Results Math., 1990, \textbf{17}, n. 3-4, 206--237.

\bibitem{Bruna} \textsc{J. Bruna}: \textit{On inverse-closed algebras of infinitely differentiable functions}, Studia Mathematica, 1980, \textbf{LXIX}, 59--68.

\bibitem{dissertationesser} \textsc{C. Esser}: \textit{Regularity of functions: {G}enericity and multifractal analysis}, PhD Thesis, Université de Liège, available online at \url{http://orbi.ulg.ac.be/bitstream/2268/174112/4/These.pdf}, 2014.
    
\bibitem{hoermander} \textsc{L. Hörmander}: The analysis of linear partial differential operators I, Distribution theory and Fourier analysis, Springer-Verlag, 2003.
    
\bibitem{Komatsu73} \textsc{H. Komatsu}: \textit{Ultradistributions. {I}. {S}tructure theorems and a characterization}, J. Fac. Sci. Univ. Tokyo Sect. IA Math., 1973, \textbf{20}, 25--105.

\bibitem{mandelbrojtbook} \textsc{S. Mandelbrojt}: Séries adhérentes, Régularisation des suites, Applications, Gauthier-Villars, Paris, 1952.

\bibitem{pietsch} \textsc{A. Pietsch}: Nuclear locally convex spaces, Ergenisse der Mathematik und ihrer Grenzgebiete, Band 66, translated from the second German edition, Springer-Verlag, 1972.

\bibitem{compositionpaper} \textsc{A. Rainer, G. Schindl}: \textit{Composition in ultradifferentiable classes}, Studia Mathematica, 2014, \textbf{224}, n. 2, 97--131.

\bibitem{characterizationstabilitypaper} \textsc{A. Rainer, G. Schindl}: \textit{Equivalence of stability properties for ultradifferentiable function classes}, accepted for publication in Rev. R. Acad. Cienc. Exactas Fis. Nat. Ser. A Math. RACSAM, available online at \url{http://arxiv.org/pdf/1407.6673.pdf}, 2014.
    
\bibitem{diploma} \textsc{G. Schindl}: Spaces of smooth functions of {D}enjoy-{C}arleman-type, Diploma Thesis, Universität Wien, available online at \url{http://othes.univie.ac.at/7715/1/2009-11-18_0304518.pdf}, 2009.

\bibitem{dissertation} \textsc{G. Schindl,}: Exponential laws for classes of {D}enjoy-{C}arleman-differentiable mappings, PhD Thesis, Universität Wien, available online at \url{http://othes.univie.ac.at/32755/1/2014-01-26_0304518.pdf}, 2014.

\bibitem{surjectivity} \textsc{J. Schmets, M. Valdivia}: \textit{On certain extension theorems in the mixed Borel setting}, J. Math. Anal. Appl., 2003, \textbf{297}, 384--403.
    
\bibitem{intersectionpaper} \textsc{J. Schmets, M. Valdivia}: \textit{Intersections of non quasi-analytic classes of ultradifferentiable functions}, Bulletin de la Société Royale des Sciences de Liège, 2008, \textbf{77}, 29--43.

\bibitem{thilliez} \textsc{V. Thilliez}: \textit{On quasi-analytic local rings}, Expo. Math., 2008, \textbf{26}, 1--23.

\end{thebibliography}
\end{document}